\documentclass[10pt]{emsart}
\usepackage{amsmath,amssymb,amstext,array}
\usepackage[latin1]{inputenc}
\usepackage{makeidx}
\makeindex

\contact[\texttt{j.m.figueroa@ed.ac.uk}]{José Figueroa-O'Farrill\\ School of
  Mathematics\\ and Maxwell Institute for Mathematical Sciences\\
  University of Edinburgh\\ United Kingdom}

\DeclareMathOperator{\Spin}{Spin}
\DeclareMathOperator{\CW}{CW}
\DeclareMathOperator{\NW}{NW}

\DeclareMathOperator{\rank}{rank}
\DeclareMathOperator{\dvol}{dvol}
\DeclareMathOperator{\im}{Im}
\DeclareMathOperator{\SU}{SU}
\DeclareMathOperator{\U}{U}
\DeclareMathOperator{\SO}{SO}

\DeclareMathOperator{\SL}{SL}
\DeclareMathOperator{\Sp}{Sp}
\DeclareMathOperator{\Mat}{Mat}
\DeclareMathOperator{\End}{End}
\DeclareMathOperator{\Hom}{Hom}
\DeclareMathOperator{\Ric}{Ric}
\DeclareMathOperator{\Riem}{Riem}

\DeclareMathOperator{\AdS}{AdS}

\DeclareMathOperator{\sS}{\mbox{\Large\textsf{S}}}
\DeclareMathOperator{\fS}{\mathfrak{S}}
\newcommand{\cyclicsum}{\sS\displaylimits}
\newcommand{\fd}{\mathfrak{d}}
\newcommand{\fg}{\mathfrak{g}}
\newcommand{\fm}{\mathfrak{m}}
\newcommand{\fn}{\mathfrak{n}}
\newcommand{\fh}{\mathfrak{h}}
\newcommand{\fk}{\mathfrak{k}}
\newcommand{\fp}{\mathfrak{p}}
\newcommand{\fso}{\mathfrak{so}}
\newcommand{\fsu}{\mathfrak{su}}
\newcommand{\fu}{\mathfrak{u}}
\renewcommand{\d}{\partial}
\newcommand{\eS}{\mathcal{S}}
\newcommand{\half}{\tfrac12}
\newcommand{\Cl}{\mathrm{C}\ell}

\newcommand{\EE}{\mathbb{E}}
\newcommand{\RR}{\mathbb{R}}
\newcommand{\HH}{\mathbb{H}}
\newcommand{\bD}{\boldsymbol{D}}
\newcommand{\ebS}{\boldsymbol{\eS}}
\newcommand{\bS}{\boldsymbol{S}}
\newcommand{\blambda}{\boldsymbol{\lambda}}

\newcommand{\1}{\boldsymbol{1}}
\newcommand{\bH}{\boldsymbol{H}}

\newcommand{\bv}{\boldsymbol{v}}
\newcommand{\bw}{\boldsymbol{w}}
\newcommand{\bx}{\boldsymbol{x}}

\newcommand{\eH}{\mathcal{H}}
\newcommand{\eL}{\mathcal{L}}
\newcommand{\eM}{\mathcal{M}}
\newcommand{\eV}{\mathcal{V}}
\newcommand{\idx}[2]{\index{#2}{\bfseries #1}}
\numberwithin{equation}{section}
\newtheorem{theorem}{Theorem}[section]
\newtheorem{corollary}[theorem]{Corollary}

\newtheorem{proposition}[theorem]{Proposition}

\theoremstyle{definition}
\newtheorem{definition}[theorem]{Definition}

\title[Lorentzian symmetric spaces in supergravity]{Lorentzian
  symmetric spaces in supergravity}

\author[José Figueroa-O'Farrill]{José Figueroa-O'Farrill}

\begin{document}

\begin{abstract}
  I will discuss the emergence of lorentzian symmetric spaces as
  supersymmetric supergravity backgrounds. I will focus on
  supergravity theories in dimension 11, 10, and 6, and will
  concentrate on the determination of the so-called maximally
  supersymmetric backgrounds, for which a classification exists up to
  local isometry.  A special class of lorentzian symmetric spaces also
  plays a rôle in the determination of parallelisable supergravity
  backgrounds in type II supergravity, which I will also summarise.\\
  (To appear in the proceedings of the programme \emph{Geometry of
    pseudo-riemannian manifolds with application to physics} hosted by
  the Erwin Schrödinger International Institute for Mathematical
  Physics.)
\end{abstract}

\begin{classification}
  Primary 53-B30; Secondary 83-50.
\end{classification}

\begin{keywords}
  Lorentzian geometry; supergravity.
\end{keywords}

\maketitle
\tableofcontents

\section{Introduction}

The purpose of this short review is to highlight the rôle of
lorentzian geometry in the supergravity limit of string theories.
That lorentzian geometry plays a rôle in such theories should not come
as a surprise, given that the supergravity theories in question
\emph{are} lorentzian theories; that is, some (if not all) of the
gravitational degrees of freedom are encoded in the form of a (local)
lorentzian metric.  What may be a little surprising is the fact that,
with notable exceptions, until relatively recently the lorentzian
nature of the solutions had not been fully exploited.  Indeed, most
early papers studying solutions to the supergravity field equations
concentrated on decomposable geometries $L \times R$, where $L$ is a
lorentzian space-form (e.g., Minkowski or (anti) de Sitter spacetime)
and $R$ a riemannian manifold, which would invariably become the focus
of the ensuing analysis.  This is not to say that such solutions are
geometrically or physically uninteresting.  In fact, they have
motivated the study of a large class of riemannian geometries which
otherwise might have remained largely in obscurity: Calabi--Yau
manifolds, manifolds with $G_2$ and $\Spin(7)$ holonomy, as well as
manifolds whose metric cones have such holonomies, e.g.,
Sasaki--Einstein and nearly Kähler manifolds, among others.  However
they do miss interesting solutions, as we will try to illustrate in
this review.

This review is organised as follows.  In Section~\ref{sec:geos} we
introduce the geometries of interest, namely lorentzian symmetric
spaces and in particular thoese which admit an absolute parallelism,
which translates into the question of which of these spaces are Lie
groups admitting bi-invariant lorentzian metrics.  In Section
\ref{sec:sugra} we discuss the geometrical aspects of supergravity
theories.  Much more could and should eventually be written about
this, but for the purposes of this review we will limit ourselves to
treat supergravity theories as collections of geometric PDEs whose
form is highly constrained, despite at first seeming ad hoc.  For
reasons explained in the body of the review, we will consider only the
following supergravity theories: eleven-dimensional supergravity,
ten-dimensional type IIB and the chiral six-dimensional $(1,0)$ and
$(2,0)$ supergravities.  In Section \ref{sec:maximal} we discuss the
classification of maximally supersymmetric solutions of the above
theories and also of ten-dimensional type IIA supergravity, which can
be obtained from eleven-dimensional supergravity via
Ka{\l}u\.za--Klein reduction, a technique we review in
Section~\ref{sec:kk}.  Finally in Section~\ref{sec:parallel} we
discuss parallelisable backgrounds in the common sector of type II
supergravity.

\subsection*{Acknowledgments}

Most of the work described in this review was obtained in
collaboration with a number of colleagues whom it is a pleasure to
thank and remember: Matthias Blau, Ali Chamseddine, Chris Hull, Kawano
Teruhiko, Patrick Meessen, George Papadopoulos, Simon Philip, Wafic
Sabra, Joan Simón, Sonia Stanciu and Yamaguchi Satoshi.  A lot of this
work was either done or started while I was a guest and/or
co-organiser of scientific programmes hosted by the Erwin Schrödinger
Institute for Mathematical Physics in Vienna and it is again my
pleasure to thank the ESI staff for all their support.  It is only
fitting that this review should appear in a volume born out of one
such programme and I would like to take this opportunity to reiterate
my thanks to Dmitri Alekseevsky and Helga Baum for the chance to
participate in it.

\section{The geometries of interest}
\label{sec:geos}

In this section we will quickly review the geometries of interest:
lorentzian symmetric spaces and, in particular, those which are
parallelisable.

\subsection{Lorentzian symmetric spaces}
\label{sec:symmetric}

We start by reviewing the classification of lorentzian symmetric
spaces.

The classification of symmetric spaces in indefinite signature is
hindered by the fact that there is no splitting theorem saying that if
the holonomy representation is reducible, the space is locally
isometric to a product.  In fact, local splitting implies both
reducibility \emph{and} a nondegeneracy condition on the factors
\cite{Wu}.  This means that one has to take into account reducible yet
indecomposable holonomy representations.  The general semi-riemannian
case is still open, but indecomposable lorentzian symmetric spaces
were classified by Cahen and Wallach \cite{CahenWallach} almost four
decades ago.  Indeed, they stated the following theorem

\begin{theorem}[Cahen--Wallach \cite{CahenWallach}]
  \label{thm:CahenWallach}
  Let $(M,g)$ be a simply-connected lorentzian symmetric space.  Then
  $M$ is isometric to the product of a simply-connected riemannian
  symmetric space and one of the following:
  \begin{itemize}
  \item $\RR$ with metric $-dt^2$;
  \item the universal cover of $n$-dimensional
    \idx{de~Sitter}{lorentzian symmetric space!de~Sitter} or
    \idx{anti~de~Sitter}{lorentzian symmetric space!anti~de~Sitter}
    spaces, where $n\geq 2$; or
  \item a \idx{Cahen--Wallach space}{lorentzian symmetric
      space!Cahen--Wallach} $\CW_n(A)$ with $n\geq 3$ and metric given
    by \eqref{eq:CWmetric} below.
  \end{itemize}
\end{theorem}
If we drop the hypothesis of simply-connectedness then this theorem
holds up to local isometry, which is the version of the theorem of
greater relevance in supergravity.

The $n$-dimensional Cahen--Wallach spaces $\CW_n(A)$ are constructed
as follows.  Let $V$ be a real vector space of dimension $n-2$ endowed
with a euclidean structure $\left<-,-\right>$.  Let $V^*$ denote its
dual. Let $Z$ be a real one-dimensional vector space and $Z^*$ its
dual.  We will identify $Z$ and $Z^*$ with $\RR$ via canonical dual
bases $\{e_+\}$ and $\{e_-\}$, respectively.  Let $A \in S^2V^*$ be a
symmetric bilinear form on $V$.  Using the euclidean structure on $V$
we can associate with $A$ an endomorphism of $V$ also denoted $A$:
\begin{equation*}
 \left< A(v), w\right> = A(v,w)\qquad\text{for all $v,w\in V$.}
\end{equation*}
We will also let $\flat:V \to V^*$ and $\sharp: V^* \to V$ denote the
musical isomorphisms associated to the euclidean structure on $V$.

Let $\fg_A$ be the Lie algebra with underlying vector space $V \oplus
V^* \oplus Z \oplus Z^*$ and with Lie brackets
\begin{equation}
  \label{eq:gA}
  \begin{aligned}[m]
    [e_-, v] &= v^\flat\\
    [e_-, \alpha] &= A(\alpha^\sharp)\\
    [\alpha, v] &= A(v,\alpha^\sharp) e_+~,
  \end{aligned}
\end{equation}
for all $v\in V$ and $\alpha \in V^*$.  All other brackets not
following from these are zero.  The Jacobi identity is satisfied by
virtue of $A$ being symmetric.  Notice that since its second derived
ideal is central, $\fg_A$ is (three-step) solvable.

Notice that $\fk_A=V^*$ is an abelian Lie subalgebra, and its
complementary subspace $\fp_A= V \oplus Z \oplus Z^*$ is acted on by
$\fk_A$.  Indeed, it follows easily from \eqref{eq:gA} that
\begin{equation*}
  [\fk_A,\fp_A] \subset \fp_A \qquad\text{and}\qquad [\fp_A,\fp_A]
  \subset \fk_A~,
\end{equation*}
whence $\fg_A = \fk_A \oplus \fp_A$ is a symmetric split.  Lastly, let
$B\in\left(S^2\fp_A^*\right)^{\fk_A}$ denote the invariant symmetric
bilinear form on $\fp_A$ defined by
\begin{equation*}
  B(v,w) = \left<v,w\right> \qquad\text{and}\qquad B(e_+,e_-)=1~,
\end{equation*}
for all $v,w\in V$.  This defines on $\fp_A$ a $\fk_A$-invariant
lorentzian scalar product of signature $(1,n-1)$.

We now have the required ingredients to construct a (lorentzian)
symmetric space.  Let $G_A$ denote the connected, simply-connected Lie
group with Lie algebra $\fg_A$ and let $K_A$ denote the Lie subgroup
corresponding to the subalgebra $\fk_A$.  The lorentzian scalar product
$B$ on $\fp_A$ induces a lorentzian metric $g$ on the space of cosets
\begin{equation*}
  M_A = G_A/K_A~,
\end{equation*}
turning it into a symmetric space.

Introducing coordinates $x^\pm, x^i$ naturally associated to
$e_\pm,e_i$, where $e_i$ is an orthonormal frame for $V$, we can write
the Cahen--Wallach metric explicitly as
\begin{equation}
  \label{eq:CWmetric}
  g = 2 dx^+ dx^- + \left(\sum_{i,j=1}^{n-2} A_{ij} x^i x^j\right)
  (dx^-)^2 + \sum_{i=1}^{n-2} \left(dx^i\right)^2~.
\end{equation}

\begin{proposition}[Cahen--Wallach \cite{CahenWallach}]
  The metric on $M_A$ defined above is indecomposable if and only if
  $A$ is nondegenerate.  Moreover, $M_A$ and $M_{A'}$ are isometric if
  and only if $A$ and $A'$ are related in the following way:
  \begin{equation*}
    A'(v,w) = c A(Ov, Ow) \qquad\text{for all $v,w\in V$,}
  \end{equation*}
  for some orthogonal transformation $O:V\to V$ and a positive scale
  $c>0$.
\end{proposition}

From this result one sees that the moduli space $\eM_n$ of
indecomposable such metrics in $n$ dimensions is given by
\begin{equation*}
  \eM_n = \left( S^{n-3} - \Delta \right) / \fS_{n-2}~,
\end{equation*}
where
\begin{equation*}
  \Delta = \left\{ (\lambda_1,\dots,\lambda_{n-2}) \in S^{n-3} \subset
    \RR^{n-2} \mid \lambda_1 \cdots \lambda_{n-2} = 0\right\}
\end{equation*}
is the singular locus consisting of eigenvalues of degenerate
$A$'s, and $\fS_{n-2}$ is the symmetric group in $n-2$ symbols, acting
by permutations on $S^{n-3} \subset \RR^{n-2}$.

\paragraph{Local isometric embeddings.}
\label{sec:embeddings}
\index{lorentzian symmetric space!isometric embedding}

Indecomposable lorentzian symmetric spaces in $d\geq 2$ are locally
isometric to algebraic varieties in pseudo-euclidean spaces.  This is
well-known for both de~Sitter and anti~de~Sitter spaces.  Indeed, let
$\kappa>0$.  Then the quadric in $\EE^{n,1}$ consisting of points
$(x_0,x_1,\dots,x_n)\in \RR^{n+1}$ such that
\begin{equation*}
  - x_0^2 + x_1^2 + x_2^2 + \cdots + x_n^2 = 1/\kappa^2
\end{equation*}
has constant sectional curvature $\kappa$ and hence is locally
isometric to a de~Sitter space, whereas the quadric in $\EE^{n-1,2}$
consisting of points $(x_0,x_1,\dots,x_n) \in \RR^{n+1}$ such that
\begin{equation*}
  -x_0^2 + x_1^2 + x_2^2 + \cdots + x_{n-1}^2 - x_n^2 =
  -1/\kappa^2
\end{equation*}
has constant section curvature $-\kappa$ and hence is locally
isometric to an anti~de~Sitter space.

Similarly, the $n$-dimensional Cahen--Wallach spaces are locally
isometric to the intersection of two quadrics in $\EE^{n,2}$.  Indeed,
let $(u^1,v^1,u^2,v^2,x^i)$, for $i=1,\dots,n-2$,
be flat coordinates in $\EE^{n,2}$ relative to which the metric takes
the form
\begin{equation}
  2 du^1 dv^1 + 2 du^2 dv^2 + \sum_{i=1}^9 dx^i dx^i~.
\end{equation}
Then the Cahen--Wallach space with matrix $A$ is locally isometric to
the induced metric on the intersection of the two quadrics
\begin{equation}
  (u^1)^2 + (u^2)^2 = 1 \qquad\text{and}\qquad 2 u^1 v^1 + 2 u^2 v^2 = 
  \sum_{i,j=1}^9 A_{ij} x^i x^j~.
\end{equation}
This was proven in \cite{Limits}.

\subsection{Lorentzian parallelisable manifolds}
\label{sec:parallelisable}

A subclass of the lorentzian symmetric spaces are the parallelisable
manifolds.

Recall that a differentiable manifold $M$ is said to admit an
\textbf{absolute parallelism} if it admits a smooth
trivialisation of the frame bundle.  Such a trivialisation consists of
a smooth global frame and hence also trivialises the tangent bundle;
whence manifolds admitting absolute parallelisms are parallelisable in
the topological sense.  The reduction theorem for connections on
principal bundles (see, for example,
\cite[Section~II.7]{KobayashiNomizu}) allows us to think of absolute
parallelisms in terms of holonomy groups of connections.  Indeed, an
absolute parallelism is equivalent to a smooth connection on the frame
bundle with trivial holonomy.  This implies, in particular, that the
connection is flat and if the manifold is simply-connected then
flatness is also sufficient.

So far these notions are purely (differential) topological and make no
mention of metrics or any other structure on the manifold.  The
question arises whether there is a metric on $M$ which is consistent
with a given absolute parallelism, so that parallel transport is an
isometry; or turning the question around, whether a given
pseudo-riemannian manifold $(M,g)$ admits an absolute parallelism
consistent with it.  In terms of connections, a consistent absolute
parallelism is equivalent to a metric connection with torsion with
trivial holonomy; or, locally, to a flat metric connection with
torsion.

Cartan and Schouten \cite{CartanSchouten1,CartanSchouten2} essentially
solved the riemannian case by generalising Clifford's parallelism on
the 3-sphere in two different ways.  The three-sphere can be
understood both as the unit-norm quaternions and also as the Lie group
$\SU(2)=\Sp(1)$.  The latter characterisation generalises to other
(semi)simple Lie groups, whereas the former gives rise to the
parallelism of the 7-sphere thought of as the unit-norm octonions.  It
follows from the results of Cartan and Schouten that a
simply-connected irreducible riemannian manifold admitting a
consistent absolute parallelism (equivalently a flat metric
connection) is isometric to one of the following: the real line, a
simple Lie group with the bi-invariant metric induced from a multiple
of the Killing form, or the round 7-sphere.

Their proofs might have had gaps which were addressed by Wolf
\cite{Wolf1,Wolf2}, who also generalised these results to arbitrary
signature, subject to an algebraic curvature condition saying that the
pseudo-riemannian manifold $(M,g)$ is of ``reductive type,'' a
condition which is automatically satisfied in the riemannian case.
(See Wolf's paper for the precise condition.)  In the case of
lorentzian signature, Cahen and Parker \cite{CahenParker} showed that
one can relax the ``reductive type'' condition; completing the
classification of absolute parallelisms consistent with a lorentzian
metric.

Wolf also showed that if one also assumes that the torsion is
parallel, then, in any signature, $(M,g)$ is locally isometric to a
Lie group with a bi-invariant metric.  In fact, as we will show below
in Section \ref{sec:flat}, one obtains the same result starting with
the weaker hypothesis that the torsion three-form is closed, which
will be the case needed in supergravity.

The results of Cahen and Parker \cite{CahenParker} actually show that
in lorentzian signature one gets for free that the torsion is
parallel.  Therefore it follows that an indecomposable lorentzian
manifold $(M,g)$ admits a consistent absolute parallelism if and only
if it is locally isometric to a lorentzian Lie group with bi-invariant
metric.\index{lorentzian symmetric space!parallelisable}

\subsection{Flat metric connections with closed torsion}
\label{sec:flat}

We will now show that a pseudo-riemannian manifold $(M,g)$ with a
flat metric connection with closed torsion three-form is locally
isometric to a Lie group admitting a bi-invariant metric.

Let $(M,g)$ be a pseudo-riemannian manifold and let $D$ be a metric
connection with torsion $T$. In other words, $D g = 0$ and for all
vector fields $X,Y$ on $M$, $T:\Lambda^2TM \to TM$ is defined by
\begin{equation*}
 T(X,Y) =  D_X Y - D_Y X - [X,Y]~.
\end{equation*}
In terms of the torsion-free Levi-Cività connection $\nabla$, we have
\begin{equation*}
  D_X Y = \nabla_X Y + \half T(X,Y)~.
\end{equation*}
Since both $Dg=0$ and $\nabla g = 0$, $T$ is skew-symmetric:
\begin{equation}
  \label{eq:skew}
  g(T(X,Y),Z) = - g(T(X,Z),Y)~,
\end{equation}
for all vector fields $X,Y,Z$ and gives rise to a \textbf{torsion
  three-form} $H\in\Omega^3(M)$, defined by
\begin{equation*}
  H(X,Y,Z) = g(T(X,Y),Z)~.
\end{equation*}
We will assume that $H$ is closed and in this section we will
characterise those manifolds for which $D$ is flat.

Let $R^D$ denote the curvature tensor of $D$, defined by
\begin{equation*}
  R^D(X,Y)Z = D_{[X,Y]}Z - D_X D_Y Z + D_Y
  D_X Z~.
\end{equation*}
Our strategy will be to consider the equation $R^D = 0$,
decompose it into types and solve the corresponding equations.  We
will find that $T$ is parallel with respect to both $\nabla$ and $D$,
and this will imply that $(M,g)$ is locally a Lie group with a
bi-invariant metric and $D$ the parallelising connection of Cartan and
Schouten \cite{CartanSchouten1}.

The curvature $R^D$ is given by
\begin{multline*}
  R^D(X,Y)Z = R(X,Y)Z - \half (\nabla_X T)(Y,Z) + \half (\nabla_Y
  T)(X,Z)\\
  - \tfrac14 T(X,T(Y,Z)) + \tfrac14 T(Y,T(X,Z))~,
\end{multline*}
where $R=R^\nabla$ is the curvature of the Levi-Cività connection.  The
tensor
\begin{equation*}
  R^D(X,Y,Z,W) := g(R^D(X,Y)Z,W)
\end{equation*}
takes the following form
\begin{multline*}
  R^D(X,Y,Z,W) = R(X,Y,Z,W)\\
  - \half g((\nabla_X T)(Y,Z),W) + \half g((\nabla_Y T)(X,Z),W)\\
  - \tfrac14 g(T(X,T(Y,Z)),W) + \tfrac14 g(T(Y,T(X,Z)),W)~,
\end{multline*}
where we have defined the Riemann tensor as usual:
\begin{equation*}
  R(X,Y,Z,W) := g(R(X,Y)Z,W)~.
\end{equation*}
Using equation \eqref{eq:skew} we can rewrite $R^D$ as
\begin{multline*}
  R^D(X,Y,Z,W) = R(X,Y,Z,W)\\
  - \half g((\nabla_X T)(Y,Z),W) + \half g((\nabla_Y T)(X,Z),W)\\
  + \tfrac14 g(T(X,W),T(Y,Z)) - \tfrac14 g(T(Y,W),T(X,Z))~,
\end{multline*}
which is manifestly skew-symmetric in $X,Y$ and in $Z,W$.  Observe
that unlike $R$, the torsion terms in $R^D$ do \emph{not} satisfy the
first Bianchi identity.  Therefore breaking $R^D$ into algebraic types
will give rise to more equations and will eventually allow us to
characterise the data $(M,g,T)$ for which $R^D = 0$.

Indeed, let $R^D = 0$ and consider the identity
\begin{equation*}
  \cyclicsum_{XYZ} R^D(X,Y,Z,W) = 0~,
\end{equation*}
where $\cyclicsum$ denotes signed permutations.  Since $R$ does obey the
Bianchi identity
\begin{equation*}
  \cyclicsum_{XYZ} R(X,Y,Z,W) = 0~,
\end{equation*}
we obtain the following identity
\begin{equation}
  \label{eq:bianchi}
  \cyclicsum_{XYZ} g((\nabla_X T)(Y,Z),W) =  - \tfrac12 \cyclicsum_{XYZ}
  g(T(W,X),T(Y,Z))~.
\end{equation}
Now we use the fact that the torsion three-form $H$ is closed, which
can be written as
\begin{multline*}
  g((\nabla_X T)(Y,Z),W) - g((\nabla_Y T)(X,Z),W)\\
  + g((\nabla_Z T)(X,Y),W) - g((\nabla_W T)(X,Y),Z)=0~,
\end{multline*}
or equivalently,
\begin{equation*}
  g((\nabla_W T)(X,Y),Z) = \half \cyclicsum_{XYZ} g((\nabla_X T)(Y,Z),W)~.
\end{equation*}
This turns equation \eqref{eq:bianchi} into
\begin{equation}
  \label{eq:bianchitoo}
  g((\nabla_W T)(X,Y),Z) = - \tfrac14 \cyclicsum_{XYZ} g(T(W,X),T(Y,Z))~.
\end{equation}
From this equation it follows that
\begin{equation*}
  g((\nabla_W T)(X,Y),Z) = - g((\nabla_X T)(W,Y),Z)~,
\end{equation*}
so that $g((\nabla_W T)(X,Y),Z)$ is totally skew-symmetric.  This
means that $\nabla H = dH = 0$, whence $H$ and hence $T$ are parallel.
Therefore equation \eqref{eq:bianchi} simplifies to
\begin{equation}
  \label{eq:Jacobi}
  \cyclicsum_{XYZ} g(T(W,X),T(Y,Z)) = 0~.
\end{equation}

Equation \eqref{eq:Jacobi} is the Jacobi identity for $T$.  Indeed,
notice that
\begin{multline*}
  g(T(W,X),T(Y,Z)) = H(W,X,T(Y,Z))\\
  = H(X,T(Y,Z),W) = g(T(X,T(Y,Z)),W)~,
\end{multline*}
whence equation \eqref{eq:Jacobi} is satisfied if and only if
\begin{equation}
  \label{eq:Jacobi}
    \cyclicsum_{XYZ} T(X,T(Y,Z)) = 0~.
\end{equation}
This means that the tangent space $T_pM$ of $M$ at every point $p$
becomes a Lie algebra where the Lie bracket is given by the
restriction of $T$ to $T_pM$.  More is true and the restriction to
$T_pM$ of the metric $g$ gives rise to an (ad-)invariant scalar
product:
\begin{equation*}
  g(T(X,Y),Z) = g(X,T(Y,Z))~.
\end{equation*}

By a theorem of Wolf \cite{Wolf1,Wolf2} (based on the earlier work of
Cartan and Schouten \cite{CartanSchouten1,CartanSchouten2}) if $(M,g)$
is complete then it is a discrete quotient of a Lie group with a
bi-invariant metric.  In general, we can say that $(M,g)$ is locally
isometric to a Lie group with a bi-invariant metric.

Indeed, since $D$ is flat, there exists locally a parallel frame
$\{\xi_i\}$ for $TM$.  Since $\xi_i$ is parallel, from the definition
of the torsion,
\begin{equation*}
  T(\xi_i,\xi_j) = - [\xi_i,\xi_j]~.
\end{equation*}
Moreover, since $T$ is parallel relative to $D$, we see that
$[\xi_i,\xi_j]$ is also parallel with respect to $D$, whence it can be
written as a linear combination of the $\xi_i$ with constant
coefficients.  In other words, they span a real Lie algebra $\fg$.
The homomorphism $\fg \to C^\infty(M,TM)$ whose image is the
subalgebra spanned by the $\{\xi_i\}$ integrates, once we choose a
point in $M$, to a local diffeomorphism $G\to M$.  This is also an
isometry if we use on $G$ the metric induced from the one on the Lie
algebra, whence we conclude that $(M,g)$ is locally isometric to a Lie
group with a bi-invariant metric.

\subsection{Bi-invariant lorentzian metrics on Lie groups}
\label{sec:lie}

In this section we will briefly review the structure of Lie groups
admitting a bi-invariant lorentzian metric.  It is well-known that
bi-invariant metrics on a Lie group are in bijective correspondence
with (ad-)invariant scalar products on the Lie algebra.  Therefore it
is enough to study those Lie algebras possessing an invariant
lorentzian scalar product.  We shall call them \idx{lorentzian Lie
  algebras}{Lie algebra!lorentzian} in this review.

It is well-known that reductive Lie algebras admit invariant scalar
products: Cartan's criterion allows us to use the Killing forms on the
simple factors and any scalar product on an abelian Lie algebra is
trivially invariant.  Another well-known example of Lie algebras
admitting an invariant scalar product are the classical doubles.  Let
$\fh$ be \emph{any} Lie algebra and let $\fh^*$ denote the dual space
on which $\fh$ acts via the coadjoint representation.  The definition
of the coadjoint representation is such that the dual pairing $\fh
\otimes \fh^* \to \RR$ is an invariant scalar product on the
semidirect product $\fh \ltimes \fh^*$ with $\fh^*$ an abelian ideal.
The Lie algebra $\fh \ltimes \fh^*$ is called the classical double of
$\fh$ and the invariant metric has split signature $(r,r)$ where $\dim
\fh = r$.

It turns out that all Lie algebras admitting an invariant scalar
product can be obtained by a mixture of these constructions.  Let
$\fg$ be a Lie algebra with an invariant scalar product
$\left<-,-\right>_\fg$.  Now let $\fh$ act on $\fg$ as skew-symmetric
derivations; that is, preserving both the Lie bracket and the scalar
product.  First of all, since $\fh$ acts on $\fg$ preserving
the  scalar product, we have a linear map
\begin{equation*}
  \fh \to \fso(\fg) \cong \Lambda^2 \fg~,
\end{equation*}
with dual map
\begin{equation*}
  c: \Lambda^2 \fg \to \fh^*~,
\end{equation*}
where we have used the invariant scalar product to identity $\fg$ and
$\fg^*$ equivariantly.  Since $\fh$ preserves the Lie bracket in
$\fg$, this map is a cocycle, whence it defines a class $[c]\in
H^2(\fg;\fh^*)$ in the second Lie algebra cohomology of $\fg$ with
coefficients in the trivial module $\fh^*$.  Let $\fg \times_c \fh^*$
denote the corresponding central extension.  The Lie bracket of $\fg
\times_c \fh^*$ is such that $\fh^*$ is central and if $X,Y\in\fg$,
then
\begin{equation*}
  [X,Y] = [X,Y]_\fg + c(X,Y)~,
\end{equation*}
where $[-,-]_\fg$ is the original Lie bracket on $\fg$.  Now $\fh$
acts naturally on $\fg \times_c \fh^*$ preserving the Lie bracket; the
action on $\fh^*$ being given by the coadjoint representation.  This
then allows us to define the \idx{double extension}{Lie algebra!double
  extension} of $\fg$ by $\fh$,
\begin{equation*}
  \fd(\fg,\fh) = \fh \ltimes (\fg \times_c \fh^*)
\end{equation*}
as a semidirect product.  Details of this construction can be found in
\cite{MedinaRevoy,FSsug}.  The remarkable fact is that $\fd(\fg,\fh)$
admits an invariant scalar product:
\begin{equation}
  \left<(X,h,\alpha), (Y,k,\beta)\right> = \left<X,Y\right>_\fg +
  \alpha(k) + \beta(h) + B(h,k)~,
\end{equation}
for all $X,Y\in\fg$, $h,k\in\fh$, $\alpha,\beta\in\fh^*$ and where $B$
is \emph{any} invariant symmetric bilinear form on $\fh$.

We say that a Lie algebra with an invariant scalar product is
\idx{indecomposable}{Lie algebra!indecomposable} if it cannot be
written as the direct product of two orthogonal ideals.  A theorem of
Medina and Revoy \cite{MedinaRevoy} (see also \cite{FSalgebra} for a
refinement) says that an indecomposable (finite-dimensional) Lie
algebra with an invariant scalar product is either one-dimensional,
simple, or a double extension $\fd(\fg,\fh)$ where $\fh$ is either
simple or one-dimensional and $\fg$ is a (possibly trivial) Lie
algebra with an invariant scalar product.  Any (finite-dimensional)
Lie algebra with an invariant scalar product is then a direct sum of
indecomposables.

If the scalar product on $\fg$ has signature $(s,t)$, then the scalar
product on the double extension $\fd(\fg,\fh)$ has signature
$(s+r,t+r)$, where $r = \dim\fh$.  This means that if we are
interested in lorentzian signature, we can double extend at most once
and by a one-dimensional $\fh$.

Therefore indecomposable lorentzian Lie algebras are either reductive
or double extensions $\fd(\fg,\fh)$ where $\fg$ has a
positive-definite invariant scalar product and $\fh$ is
one-dimensional.  In the reductive case, indecomposability means that
it has to be simple, whereas in the latter case, since the scalar
product on $\fg$ is positive-definite, $\fg$ must be reductive.  A
result of \cite{FSsug} (see also \cite{FSalgebra}) then says that any
semisimple factor in $\fg$ splits off resulting in a decomposable Lie
algebra.  Thus if the double extension is to be indecomposable, $\fg$
must be abelian.  In summary, an indecomposable lorentzian Lie algebra
is either simple or a double extension of an abelian Lie algebra by a
one-dimensional Lie algebra and hence solvable (see, e.g.,
\cite{MedinaRevoy}).

In summary, an indecomposable lorentzian Lie algebra is either
isomorphic to $\fso(1,2)$ with (a multiple of) the Killing form, or
else is solvable and can be described as a double extension
$\fd_{2n+2}:= \fd(\EE^{2n},\RR)$ of the abelian Lie algebra $\EE^{2n}$
with the (trivially invariant) euclidean ``dot'' product by a
one-dimensional Lie algebra acting on $\EE^{2n}$ via a non-degenerate
skew-symmetric linear map $J : \EE^{2n} \to \EE^{2n}$.  Let $\omega
\in \Lambda^2(\EE^{2n})^*$ denote the associated $2$-form:
$\omega(\bv,\bw) = \left<\bv, J\bw\right>$.

More concretely, the double extension $\fd_{2n+2}$ has underlying
vector space $V = \EE^{2(d-1)} \oplus \RR \oplus \RR$, and if
$(\bv,v^-,v^+),(\bw, w^-, w^+) \in V$, then their Lie bracket is
given by
\begin{equation*}
  [(\bv,v^-,v^+), (\bw, w^-, w^+)] = (v^- J(\bw) - w^- J(\bv), 0, \bv
  \cdot J(\bw))
\end{equation*}
and their inner product follows by polarisation from
\begin{equation*}
  \left|(\bv,v^-,v^+)\right|^2 = \bv \cdot \bv + 2 v^+ v^- + b
  (v^-)^2~,
\end{equation*}
where $b\in \RR$ is arbitrary.  One can however always set $b=0$ via a
Lie algebra automorphism and we will do so here; although there are
situations when one may wish to retain this freedom.

The unique simply-connected Lie group with Lie algebra
$\fd_{2n+2}$ is a solvable ($2n+2$)-dimensional Lie group admitting a
bi-invariant metric
\begin{equation}
  \label{eq:cwmetric}
  ds^2 = 2 dx^+ dx^- - \left<J\bx,J\bx\right> (dx^-)^2 +
  \left<d\bx,d\bx\right>~,
\end{equation}
relative to natural coordinates $(\bx, x^-, x^+)$.  The parallelising
torsion has $3$-form
\begin{equation}
  \label{eq:partor}
  H = dx^- \wedge \omega~.
\end{equation}

The non-degenerate skew-symmetric endomorphism $J$ can be brought to a
Jordan normal form consisting of nonzero $2\times 2$ blocks via an
orthogonal transformation.  The skew-eigenvalues
$\lambda_1,\dots,\lambda_n$, which are different from zero, can be
arranged so that they obey: $0 < \lambda_1 \leq \lambda_2 \leq \cdots
\leq \lambda_n$.  Finally a positive rescaling of $J$ can be absorbed
into reciprocal rescalings of $x^\pm$, so that we can set $\lambda_n$,
say, equal to $1$ without loss of generality.  Therefore we see that
the moduli space of metrics \eqref{eq:cwmetric} is given by an
($n-1$)-tuple $\blambda = (\lambda_1,\ldots,\lambda_{n-1})$ where $0 <
\lambda_1 \leq \cdots \leq \lambda_{n-1} \leq 1$.  It is clear that
they are particular cases of the Cahen--Wallach spaces discussed in
Section~\ref{sec:symmetric}.

\section{Supergravity}
\label{sec:sugra}

Supergravity is one of the later jewels of 20th century theoretical
physics.  It started out as an attempt to `gauge' the supersymmetry of
certain quantum field theories, but it was quickly realised that it
provides a nontrivial extension of Einstein gravity.  Supergravity
theories are fairly rigid---their structure dictated largely by the
representation theory of the spin groups.  A good modern review of the
structure of supergravity theories is \cite{ToineReview}.  It is fair
to say, however, that supergravity theories are still somewhat
mysterious to most mathematicians and much remains to be done to make
this beautiful chapter of modern mathematical physics accessible to a
larger mathematical audience.  That, however, is a task for a
different occasion.  For our present purposes, each supergravity
theory will be a collection of geometric PDEs and our interest will be
in finding special types of solutions.

We shall be interested uniquely in lorentzian supergravity theories in
dimension $d\geq 4$.  There are supergravity theories in lower
dimensions and in other metric signatures, but we will not discuss
them here.  Neither will we discuss other types of supergravity
theories: heterotic, gauged, conformal, massive,...  The two-volume
set \cite{SalamSezgin} reprints many of the foundational supergravity
papers.

For reasons which are well-known, namely the otherwise non-existence
of nontrivial interacting theories, the dimension of the spacetime
will be bounded above by $11$.  Apart from the dimension of the
spacetime, the other important invariant is the ``number of
supercharges'', denoted $n$, which is an integer multiple of the
dimension of the smallest irreducible \emph{real} spinor
representation in that spacetime dimension.  For dimension $\geq 4$
the number of supercharges ranges from $4$ to $32$.

In Table \ref{tab:supergravities}, which is borrowed from
\cite{ToineReview}, we tabulate the different supergravity theories in
$d\geq 4$.  The seemingly baroque notation is not too important: M
refers to the unique eleven-dimensional supergravity theory
\cite{Nahm,CJS} which is a low-energy limit of M-theory (hence the
name), types I \cite{NilssonN=1,ChamseddineN=1}, IIA
\cite{CampbellWestIIA,HuqNamazieIIA,GianiPerniciIIA} and IIB
\cite{SchwarzIIB,SchwarzWestIIB,HoweWestIIB} supergravities are the
low-energy limits of the similarly named string theories, whereas the
notation $N=n$ or $(p,q)$ is historical and denotes the multiplicity
of the spinor (or half-spinor) representations in the corresponding
supersymmetry algebra.  The original supergravity theory
\cite{FvNF,DeserZumino} is the four-dimensional $N=1$ theory.  The top
entry in each column has been highlighted to indicate that upon
dimensional reduction it gives rise to all the theories below it in
the same column.  As we will explain below, this means that for many
purposes, especially the classification of solutions, it is generally
enough to understand the `top' theories and, indeed, we will
concentrate on those.

\begin{table}[h!]
  \centering
  \begin{tabular}{|*{8}{>{$}c<{$}|} }
    \hline
    d\downarrow~n\to & \multicolumn{4}{c|}{32} & \multicolumn{2}{c|}{24}  &
    20 \\
    \hline
    11  & \textbf{M} & \multicolumn{3}{c|}{} &\multicolumn{2}{c|}{} & \\
    10  & \text{IIA} & \textbf{IIB} & \multicolumn{2}{c|}{}&\multicolumn{2}{c|}{} &
    \\
    9  &  \multicolumn{2}{c|}{N=2} &\multicolumn{2}{c|}{ }&
    \multicolumn{2}{c|}{ }  & \\
    8  &  \multicolumn{2}{c|}{N=2}&\multicolumn{2}{c|}{ }&
    \multicolumn{2}{c|}{ }  & \\
    7  &  \multicolumn{2}{c|}{N=4} &\multicolumn{2}{c|}{
    }&\multicolumn{2}{c|}{ }  & \\
    6  &
    \multicolumn{2}{c|}{(2,2)}&\boldsymbol{(3,1)}&\boldsymbol{(4,0)}
    &\boldsymbol{(2,1)} & \boldsymbol{(3,0)}& \\
    5  &  \multicolumn{4}{c|}{N=8}  &\multicolumn{2}{c|}{N=6}  & \\
    4  &  \multicolumn{4}{c|}{N=8}  & \multicolumn{2}{c|}{N=6}&
    \boldsymbol{N=5} \\
    \hline
  \end{tabular}\\[20pt]
  \begin{tabular}{|*{6}{>{$}c<{$}|} }
    \hline
    d\downarrow~n\to & \multicolumn{2}{c|}{16}  & 12 & 8 & 4  \\
    \hline
    10  & \textbf{I} &  &  &  &  \\
    9  & N=1  &  &  &  &  \\
    8  & N=1 &  &  &  &  \\
    7  & N=2 & &  &  &  \\
    6  & (1,1) &  \boldsymbol{(2,0)} &  & \boldsymbol{(1,0)} &  \\
    5  & \multicolumn{2}{c|}{N=4}&  & N=2 &  \\
    4  & \multicolumn{2}{c|}{N=4}  & \boldsymbol{N=3} & N=2 & \boldsymbol{N=1}\\
    \hline
  \end{tabular}
  \label{tab:supergravities}
  \caption{Lorentzian supergravity theories in $d\geq 4$}
\end{table}

Indeed, supergravity theories in different dimensions may be related by a
procedure known as Ka{\l}u\.za--Klein reduction.  This can be read off
from Table \ref{tab:supergravities}: any supergravity theory in the
table can be obtained by Ka{\l}u\.za--Klein reduction from any theory
sitting above it in the same column.  In practice this means that a
solution to any of the supergravity theories in the table can be
lifted to a solution of any theory above it in the same column, should
there be any.  Conversely, any solution of a supergravity theory which
is invariant under a one-dimensional Lie group gives rise to a local
solution (and indeed global if the action is free and proper) of the
supergravity theory immediately below it in the same column.  We shall
be particularly interested in the reduction from $d=11$ supergravity
to $d=10$ type IIA supergravity, and in the reduction and subsequent
truncation from the $d=6$ $(1,0)$ supergravity to minimal $d=5$ $N=2$
supergravity.

We shall be interested in solutions of the field equations coming from
these supergravity theories.  Such a solution is described in
geometric terms by the following data:
\begin{itemize}
\item a $d$-dimensional lorentzian spin manifold $(M,g)$ with a
  (possibly twisted) real rank $n$ spinor bundle $\eS \to M$, and
\item certain additional geometric data, which will be different in
  each supergravity theory, consisting of differential forms or, more
  generally, sections of certain fibre bundles over $M$,
\end{itemize}
all subject to field equations which generalise the coupled
Einstein--Maxwell equations familiar from four-dimensional Physics.

The above geometric data defines a connection $D$ on the spinor bundle
$\eS$ as well as a (possibly empty) set of endomorphisms of $\eS$.
Together they define a class of sections of $\eS$, parallel with
respect to $D$ and in the kernel of the endomorphisms, which are
called Killing spinors.

We will be particularly interested in the cases where the connection
$D$ is flat, so that it admits the maximum number of parallel
sections.  In this case, the field equations are automatically
satisfied.  In general, the field equations are intimately related to
integrability conditions for the existence of parallel sections of
$D$.

\subsection{Eleven-dimensional supergravity}
\label{sec:d=11}

Eleven-dimensional supergravity was predicted by Nahm \cite{Nahm} and
constructed soon thereafter by Cremmer, Julia and Scherk \cite{CJS}.
We will only be concerned with the bosonic equations of motion.  The
geometrical data consists of $(M, g, F)$ where $(M,g)$ is an
eleven-dimensional lorentzian manifold with a spin structure and $F
\in \Omega^4(M)$ is a closed $4$-form.  The equations of motion
generalise the Einstein--Maxwell equations in four dimensions.  The
Einstein equation relates the Ricci curvature to the energy momentum
tensor of $F$.  More precisely, the equation is
\begin{equation}
  \label{eq:einstein}
    \Ric(g) = T(g,F)
\end{equation}
where the symmetric tensor
\begin{equation*}
  T(X,Y) = \half \left< \imath_X F, \imath_Y F\right> - \tfrac16
  g(X,Y) |F|^2~,
\end{equation*}
is related to the energy-momentum tensor of the (generalised) Maxwell
field $F$.  In the above formula, $\left<-,-\right>$ denotes the
scalar product on forms, which depends on $g$, and $|F|^2 =
\left<F,F\right>$ is the associated (indefinite) norm.  The
generalised Maxwell equations are now nonlinear:
\begin{equation}
  \label{eq:maxwell}
   d \star F = - \half F \wedge F~.
\end{equation}

\begin{definition}
  A triple $(M,g,F)$ satisfying the equations \eqref{eq:einstein} and
  \eqref{eq:maxwell} is called a (bosonic)
  \idx{background}{supergravity!$d{=}11$!background} of
  eleven-dimensional supergravity.
\end{definition}

Let $\eS$ denote the bundle of spinors on $M$.  It is a real vector
bundle of rank $32$ with a spin-invariant symplectic form
$\left(-,-\right)$.  A differential form on $M$ gives rise to an
endomorphism of the spinor bundle via the composition
\begin{equation*}
  c: \Lambda T^*M  \xrightarrow{\cong} \Cl(T^*M) \to \End\eS~,
\end{equation*}
where the first map is the bundle isomorphism induced by the vector
space isomorphism between the exterior and Clifford algebras, and the
second map is induced from the action of the Clifford algebra
$\Cl(1,10)$ on the spinor representation $S$ of $\Spin(1,10)$.  In
signature $(1,10)$ one has the algebra isomorphism
\begin{equation*}
  \Cl(1,10) \cong \Mat(32,\RR) \oplus \Mat(32,\RR)~,
\end{equation*}
hence the map $\Cl(1,10) \to \End S$ has kernel.  In other words, the
map $c$ defined above involves a choice.  This comes down to choosing
whether the (normalised) volume element in $\Cl(1,10)$ acts as $\pm$
the identity.  In our conventions, the volume element acts as minus
the identity.

\begin{definition}
  We say that a background $(M,g,F)$ is
  \idx{supersymmetric}{supergravity!$d{=}11$!supersymmetric
    background} if there exists a nonzero spinor $\varepsilon \in
  \Gamma(\eS)$ which is parallel with respect to the supersymmetric
  connection
  \begin{equation*}
    D : \Gamma(\eS) \to \Gamma( T^*M \otimes \eS)
  \end{equation*}
  defined, for all vector fields $X$, by
  \begin{equation}
    \label{eq:connd=11}
    D_X \varepsilon  = \nabla_X \varepsilon + \Omega_X(F)
    \varepsilon~,
  \end{equation}
  where $\nabla$ is the spin connection and $\Omega(F):TM \to \End\eS$
  is defined by
  \begin{equation*}
    \Omega_X(F) = \tfrac1{12} c(X^\flat \wedge F) - \tfrac16
    c(\imath_X F)~,
  \end{equation*}
  with $X^\flat$ the one-form dual to $X$.
\end{definition}

A nonzero spinor $\varepsilon$ which is parallel with respect to $D$
is called a \idx{Killing spinor}{supergravity!$d{=}11$!Killing
  spinor}.  This is a generalisation of the usual geometrical notion
of Killing spinor (see, for example, \cite{BFGK}).  The name is apt
because Killing spinors are ``square roots'' of Killing vectors.
Indeed, one has the following

\begin{proposition}
  Let $\varepsilon_i$, $i=1,2$ be Killing spinors: $D \varepsilon_i =
  0$.  Then the vector field $V$ defined, for all vector fields $X$,
  by
  \begin{equation*}
    g(V,X) = \left( \varepsilon_1, X \cdot \varepsilon_2\right)
  \end{equation*}
  is a Killing vector and moreover \cite{GauPak} preserves $F$.
\end{proposition}

The fundamental object in eleven-dimensional supergravity is the
connection $D$, whose curvature encodes the field equations.  Indeed,
the field equations are equivalent \cite{GauPak,Preons} to
\begin{equation*}
  e^i \cdot R^D_{X,e_i} = 0\qquad\text{for every vector field $X$,}
\end{equation*}
where $(e_i)$ is an orthonormal frame and $(e^i)$ is dual coframe and
$\cdot$ is Clifford multiplication.

Alas, $D$ is not induced from a connection on the tangent bundle and
in fact, it does not even preserve the symplectic structure.
Nevertheless one has the following

\begin{proposition}\cite{HullHolonomy}
  The holonomy of $D$ is contained in $\SL(32,\RR)$.
\end{proposition}

An important open problem is to determine the possible holonomy groups
of $D$ subject to the field equations.  In a way, the field equations
play the rôle of the torsion-free condition in the holonomy problem
for affine connections.  Except for the above result there are no
other results of a general nature and although the infinitesimal
holonomy of a number of solutions are known \cite{DuffLiuHol,BDLWHol},
a general pattern has yet to emerge.

A coarser invariant than the holonomy of $D$ is the dimension of its
kernel; that is, the dimension of the space of Killing spinors. It is
customary to write this as a fraction
\begin{equation*}
  \nu = \frac{\dim \{\text{Killing spinors}\}}{\rank \eS}
\end{equation*}
which in this case is of the form $k/32$ for some integer
$k=0,1,\dots,32$.

In Section~\ref{sec:maxd=11} we will review the classification of
those bacgrounds with $\nu = 1$; that is, those backgrounds where $D$
is flat.  We will see that they are all given by lorentzian symmetric
spaces.  In fact, it was shown in \cite{FMPHom} that if $\nu > \tfrac34$,
then $M$ is locally homogeneous and moreover it was conjectured that
there exist backgrounds with $\nu = \tfrac34$ which are not locally
homogeneous; although at present none have been constructed.  At the
other end of the spectrum, the general form of $(g,F)$ which admit
(at least) one Killing spinor is known \cite{GauPak,GauGutPak}.

\subsection{Ten-dimensional IIB supergravity}
\label{sec:IIB}

Ten-dimensional IIB supergravity
\cite{SchwarzIIB,SchwarzWestIIB,HoweWestIIB} is somewhat more
complicated than eleven-dimensional supergravity due to the
proliferation of dynamical fields and the fact that it cannot be
obtained by Ka{\l}u\.za--Klein reduction from any higher-dimensional
supergravity theory.

A type IIB supergravity background is described by the geometric data
we describe presently.  First of all, we have a ten-dimensional
lorentzian spin manifold $(M,g)$ together with a self-dual $5$-form
$F$.  Now let $\eH$ be the complex upper half-plane, thought of as the
riemannian symmetric space $\SU(1,1)/\U(1)$ and let $\tau: M\to \eH$
be a smooth map.  We may think of $\SU(1,1)$ as the total space of a
principle circle bundle over $\eH$ and we let $\eL\to\eH$ denote the
associated complex line bundle.  Let $L_\tau = \tau^*\eL$ denote the
pull-back bundle over $M$.  Choosing a section $\sigma : \eH \to
\SU(1,1)$, we may pull back to $M$ the left-invariant Maurer--Cartan
form on $\SU(1,1)$: its component along $\fu(1)$ defines a connection
$A$ on $L_\tau$, whereas the component perpendicular to $\fu(1)$,
relative to the invariant lorentzian scalar product on $\fsu(1,1)$,
defines a one-form $B$ on $M$ with values in $L_\tau^2$.  Both $A$ and
$B$ can be written explicitly in terms of $\tau$.  Indeed, if we let
$z = (\tau -i)/(\tau + i)$ be the Cayley transform of $\tau$, so that
$|z|<1$, then there is a choice of section $\sigma$ for which
\begin{equation*}
  A = \frac{\im(z d\overline z)}{1-|z|^2} \qquad\text{and}\qquad
  B = \frac{dz}{1-|z|^2}~.
\end{equation*}
Finally, let $G$ be a $3$-form on $M$ with values in $L_\tau$.  On the
bundles $\Lambda^p T^*M \otimes L_\tau^q$ we have connections
$\nabla^{p,q}$ obtained from the Levi-Cività connection on $TM$ (and
hence the tensor bundles) and the connection $A$ on $L_\tau$ (and
hence its powers).  We will let
\begin{equation*}
  d^{\nabla^{p,q}} : \Omega^p(M;L_\tau^q) \to \Omega^{p+1}(M;L_\tau^q)
  \quad\text{and}\quad
  \delta^{\nabla^{p,q}} : \Omega^p(M;L_\tau^q) \to
  \Omega^{p-1}(M;L_\tau^q)
\end{equation*}
denote the associated differential and co-differential on
$L_\tau^q$-valued differential forms.  We will let $\left<-,-\right>$
denote the natural pairing
\begin{equation*}
  \Omega^p(M;L_\tau^q) \otimes \Omega^p(M;L_\tau^{q'}) \to
  \Omega^0(M;L_\tau^{q+q'})
\end{equation*}
induced from the metric $g$.  With these notational remarks behind us,
we can finally define a IIB supergravity background.

\begin{definition}
  The data $(M,g,\tau,F,G)$ described above defines a IIB supergravity
  \idx{background}{supergravity!$d{=}10$ IIB!background} provided that
  the following equations are satisfied:
  \begin{align*}
    \delta^{\nabla^{1,2}} B &= \tfrac14 |G|^2 \\
    \delta^{\nabla^{3,1}} G(X,Y) &= \left<B,\imath_X\imath_Y \overline
      G \right> - \tfrac{2i}3 \left<\imath_X\imath_YF,G\right> \\
    d^{\nabla^{3,1}} G &= - B \wedge \overline G \\
    d^{\nabla^{5,0}} F &= \tfrac{i}{8} G \wedge \overline G \\
    \Ric(X,Y) &= B(X)\overline B(Y) + B(Y)\overline B(X) + 4
    \left<\imath_XF,\imath_YF\right>\\
    & \quad {} + \tfrac14 \left(\left<\imath_XG,\imath_Y\overline G\right> +
      \left<\imath_YG,\imath_X\overline G\right>\right) -
    \tfrac18 \left<G,\overline G\right> g(X,Y)~.
  \end{align*}
\end{definition}

Let $\eS_\pm$ denote the half-spinor bundles over $M$.  They are real,
symplectic and have rank $16$.  Let $\eS := \eS_- \otimes
L_\tau^{1/2}$, where $L_\tau^{1/2}$ is the square-root bundle of
$L_\tau$.  Let $\overline\eS := \eS_- \otimes L_\tau^{-1/2}$.  Notice
that if $\varepsilon \in \Omega^0(M;\eS)$, then $\overline\varepsilon
\in \Omega^0(\overline\eS)$.  Furthermore the Clifford action of
differential forms on spinors extends to an action
\begin{equation*}
  c : \Omega^p(M;L_\tau^q) \to \Hom\left(\Omega^0\left(M;\eS_\pm \otimes
  L_\tau^r\right), \Omega^0\left(M;\eS_{(-1)^p} \otimes
  L_\tau^{r+q}\right)\right)~.
\end{equation*}
Similarly we have a connection $\nabla^s$ acting on
$\Omega^0(M;\eS_\pm \otimes L_\tau^s)$ which is defined using the spin
connection and the connection $A$ on $L_\tau$.  We are now in a
position to define a type IIB supergravity Killing spinor.

\begin{definition}
  A IIB supergravity \idx{Killing spinor}{supergravity!$d{=}10$ IIB!Killing
    spinor} is a nonzero section $\varepsilon$ of $\eS$ satisfying the
  following two conditions
  \begin{align}
    c(B)\overline \varepsilon &= \tfrac14 c(G)\varepsilon \label{eq:dilatinoIIB}\\
    \nabla^{1/2}_X \varepsilon &= - \tfrac{i}4 c(F) c(X^\flat)
    \varepsilon - \tfrac1{32} \left( c(\imath_X G) - 2 c(X^\flat
      \wedge G)\right) \overline \varepsilon~.\label{eq:gravitinoIIB}
  \end{align}
\end{definition}

Just like in eleven-dimensional supergravity, IIB Killing spinors are
square roots of Killing vectors.  Indeed, the image of the natural map
\begin{equation*}
  \Omega^0(M;\eS) \otimes \Omega^0(M;\overline\eS) \to \text{vector
    fields}
\end{equation*}
consists of Killing vectors which in addition preserve the geometric
data of a background.  Again it is possible to show that if the space
of Killing spinors of a supersymmetric background of type IIB
supergravity has (real) dimension $>24$, then the background is
locally homogeneous \cite{EHJGMHom}.

Type IIB supergravity backgrounds are acted upon by $\SU(1,1)$, which
is the \idx{duality group}{supergravity!$d{=}10$ IIB!duality group} of type IIB
supergravity.  The metric $g$ and the five-form $F$ are
$\SU(1,1)$-invariant, whereas $\SU(1,1)$ acts on $z$ (hence on $\tau$)
via fractional linear transformations:
\begin{equation*}
  \begin{pmatrix}
    a & b \\ \overline b & \overline a
  \end{pmatrix} \cdot z = \frac{a z + b}{\overline b z + \overline
    a}~.
\end{equation*}
Moreover, the bundle $\eL \to \eH$ is a homogeneous bundle of
$\SU(1,1)$ hence there is an action of $\SU(1,1)$ on sections of $\eL$
and its powers.  Putting these two actions together we see that
$\gamma\in\SU(1,1)$ sends sections of $L_\tau^p$ (and also differential
forms and spinors with values in such a bundle) to sections of
$L_{\gamma\tau}^p$.

The classification of the maximally supersymmetric background will be
presented in Section~\ref{sec:maxIIB}, which is based on the papers
\cite{FOPMax} and \cite{FOPPluecker}.  In the opposite extreme, there
has been steady progress recently on the determination of the general
form of the backgrounds admitting some supersymmetry
\cite{GGP1,GGP2,GGPR}.

\subsection{Six-dimensional $(2,0)$ and $(1,0)$ supergravities}
\label{sec:sixdim}

We start by describing the field content and Killing spinor equations
of $(1,0)$ \cite{NishinoSezgin10} and $(2,0)$
\cite{Townsend20,Tanii20} chiral supergravities in six dimensions.  We
start as usual by describing the relevant spinorial representations.
The spin group $\Spin(1,5) \cong \SL(2,\HH)$, whence the irreducible
spinorial representations are quaternionic of complex dimension $4$.
There are two inequivalent representations $S_\pm$ which are
distinguished by their chirality. Let $S_1$ denote the fundamental
representation of $\Sp(1)$: it is a quaternionic representation of
complex dimension $2$, and similarly let $S_2$ denote the fundamental
representation of $\Sp(2)$, which is a quaternionic representation of
complex dimension $4$.  The tensor products $S_+ \otimes S_1$ and $S_+
\otimes S_2$ are complex representations of $\Spin(1,5) \times \Sp(1)$
and $\Spin(1,5) \times \Sp(2)$, respectively, with a real structure.
We will let
\begin{equation*}
  S = [S_+ \otimes S_1] \qquad\text{and}\qquad \bS = [S_+ \otimes
  S_2]
\end{equation*}
denote the underlying real representations.  Clearly $S$ is a real
representation of dimension $8$ and $\bS$ is a real representation of
dimension $16$.  If $(M,g)$ is a six-dimensional lorentzian spin
manifold, then we will let $\eS$ and $\ebS$ denote the bundles of
spinors associated with the representations $S$ and $\bS$,
respectively.  The groups $\Sp(1)$ and $\Sp(2)$ are the
\idx{R-symmetry}{supergravity!$d{=}6$!R-symmetry groups} groups of
these supergravity theories.

\begin{definition}
  A $(1,0)$ supergravity \idx{background}{supergravity!$d{=}6$!$(1,0)$
    background} consist of a six-dimensional lorentzian spin manifold
  $(M,g)$ together with a closed antiself-dual $3$-form $H$ subject to
  the Einstein equation
  \begin{equation*}
    \Ric(X,Y) = -\tfrac14 \left<\imath_XH,\imath_YH\right>~.
  \end{equation*}
  Such a background is said to be
  \idx{supersymmetric}{supergravity!$d{=}6$!$(1,0)$ supersymmetric
    background} if there are nonzero sections $\varepsilon$ of $\eS$
  obeying
  \begin{equation}
    \label{eq:conn10}
    D_X \varepsilon := \nabla_X \varepsilon - \half c(\imath_X H)
    \varepsilon = 0~,
  \end{equation}
  for all vector fields $X$, where $c: \Omega(M) \to \Cl(T^*M) \to
  \End(S)$ is the action of forms on sections of $S$, and $\nabla$ is
  induced from the Levi-Cività connection.
\end{definition}

We remark that the connection $D$ in equation \eqref{eq:conn10} is
induced from a spin connection with torsion three-form $H$.

Similarly for $(2,0)$ supergravity, we have the following

\begin{definition}
  A (2,0) supergravity \idx{background}{supergravity!$d{=}6$!$(2,0)$
    background} consists of a
  six-dimensional lorentzian spin manifold $(M,g)$, a $V$-valued
  closed antiself-dual $3$-form $\bH$, where $V$ is the
  five-dimensional real representation of the R-symmetry group
  $\Sp(2)\cong \Spin(5)$ together with a $\Sp(2)$-invariant scalar
  product, subject to the Einstein equation
  \begin{equation*}
    \Ric(X,Y) = -\tfrac14 \left<\imath_X\bH,\imath_Y\bH\right>~,
  \end{equation*}
  where $\left<-,-\right>$ now also includes the $\Sp(2)$-invariant
  inner product on $V$.  Such a background is said to be
  \idx{supersymmetric}{supergravity!$d{=}6$!$(2,0)$ supersymmetric
    background} if there are nonzero sections $\varepsilon$ of $\ebS$
  obeying
  \begin{equation}
    \label{eq:conn20}
    \bD_X \varepsilon = \nabla_X \varepsilon - \half c(\imath_X \bH)
    \varepsilon = 0~,
  \end{equation}
  for all vector fields $X$ and where $c: \Omega(M;V) \to \Cl(T^*M)
  \otimes \Cl(V) \to \End(\ebS)$ is the action of $V$-valued forms on
  sections of $\ebS$.
\end{definition}

Notice that in $(2,0)$ supergravity, the anti-selfduality of $\bH$
imply that $\bH \wedge \bH = 0$ in $\Omega^6(M; \Lambda^2 V)$.

Maximal supersymmetry implies that the connections $D$ acting on $\eS$
and $\bD$ on $\ebS$ are flat.  In the case of $(1,0)$ supergravity,
$D$ is a spin connection with torsion and maximally supersymmetric
solutions correspond to six-dimensional lorentzian manifolds admitting
a flat metric connection with anti-selfdual closed torsion three-form.
We saw in Section \ref{sec:flat} that $(M,g)$ is locally isometric to
a Lie group with a bi-invariant lorentzian metric.  In the case of
$(2,0)$ supergravity, $\bD$ does not have such an obvious geometrical
interpretation, but it is proven in \cite{CFOSchiral} that, up to the
natural action of the R-symmetry group, the $(2,0)$ maximally
supersymmetric backgrounds are in one-to-one correspondence with those
of $(1,0)$ supergravity.

The general form of a supersymmetric background in $(1,0)$
supergravity has been obtained in \cite{GMR}, who in particular also
determine the maximally supersymmetric backgrounds by a different
method, closely related to the one in \cite{FOPMax,FOPPluecker}.

\section{Maximally supersymmetric backgrounds}
\label{sec:maximal}

In this section we review the known results about maximally
supersymmetric backgrounds in a number of the more interesting
supergravity theories.  Several of the theories under the
consideration will be tackled directly: $d=11$ supergravity,
$d=10$ IIB supergravity and the $d=6$ supergravities, whereas the
maximally supersymmetric backgrounds of $d=10$ IIA and $d=5$ $N=2$
supergravities will be obtained from those $d=11$ and $d=6$
supergravities by the technique of Ka{\l}u\.za--Klein reduction.  As
this technique is very useful, we will review it briefly now.

\subsection{Ka{\l}u\.za--Klein reduction}
\label{sec:kk}

In this section we will briefly review the geometric underpinning of
Ka{\l}u\.za--Klein reduction.  We start with a supergravity background
$(M,g,F,\dots)$ which is invariant under a one-dimensional Lie group
$\Gamma$, acting freely and properly on $M$ by isometries which in
addition preserve any other geometric data $F,\dots$.  We shall let
$\xi$ denote a Killing vector field for the $\Gamma$-action.  Since
the action is free, $\xi$ is nowhere vanishing.  We will also assume
that $\xi$ is spacelike; although this is not strictly necessary and
indeed time-like reductions can be quite useful, especially in the
context of topological field theories.

The original spacetime $M$ is to be thought of as the total space of a
principal $\Gamma$-bundle $\pi:M \to N = M/\Gamma$, where $\pi$ the
map taking a point in $M$ to the $\Gamma$-orbit on which it lies.  At
every point $p$ in $M$, the tangent space $T_pM$ of $M$ at $p$
decomposes into two orthogonal subspaces: $T_p M = \eV_p \oplus
\eH_p$, where the \textbf{vertical subspace} $\eV_p = \ker \pi_*$
consists of those vectors tangent to the $\Gamma$-orbit through $p$,
and the \textbf{horizontal subspace} $\eH_p=\eV_p^\perp$ is its
orthogonal complement relative to the metric $g$.  The resulting
decomposition is indeed a direct sum by virtue of the
nowhere-vanishing of the norm of $\xi$, whose value at $p$ spans
$\eV_p$ for all $p$.  The derivative map $\pi_*$ sets up an
isomorphism between $T_p M$ and $T_q N$, where $\pi(p) = q$.  As is
well-known, there is a unique metric on $N$ for which this isomorphism
is also an isometry and for which the map $\pi$ is a riemannian
submersion.  We will call this metric $h$.

The horizontal sub-bundle $\eH$ gives rise to a connection one-form
$\alpha$ on $M$ such that $\eH = \ker \alpha$ and such that
$\alpha(\xi) = 1$.  We remark that $\alpha$ is invariant, so that
$\eL_\xi\alpha = 0$.  This means that the curvature $2$-form $d\alpha$
is both invariant and \textbf{horizontal}---that is, $\imath_\xi d\alpha
= 0$.  Such forms are called \textbf{basic} and it is a basic fact that
they define forms on $N$.  Hence $d\alpha$ defines a $2$-form on $N$.

Finally the norm $|\xi|$ of the Killing vector is itself
$\Gamma$-invariant and hence defines a function on $N$.  Since $\xi$
is spacelike, this function is positive and hence it is convenient to
write it as the exponential of a real valued function $\phi:N \to \RR$
which is (up to a constant multiple) called the \textbf{dilaton}.

In summary, and omitting the pull-backs on $h$ and $\phi$, we can
write the metric $g$ as
\begin{equation*}
  g = h + e^{2\phi} \alpha^2~.
\end{equation*}

The other geometric data also reduces.  For example, if $F$ is an
invariant differential form on $M$, it gives rise to two differential
forms on $N$ simply by decomposing
\begin{equation*}
  F = G - \alpha \wedge H~,
\end{equation*}
where $\alpha$ is the connection one-form defined above.  The forms
$G$ and $H$ are basic and hence define differential forms on $N$.
Indeed, it is clear from the above expression that $H = -\imath_\xi
F$, so that it is manifestly horizontal.  Invariance of $F$ means that
$H$ is closed, whence it is also invariant.  Finally, we observe that
$G$ is also basic.  It is manifestly horizontal, and invariance
follows by a simple calculation using that $G$, $H$ and $d\alpha$ are
horizontal.

\subsection{Eleven-dimensional supergravity}
\label{sec:maxd=11}

Maximal supersymmetry implies the flatness of the supersymmetric
connection (\ref{eq:connd=11}).  Calculating the curvature of this
connection and separating into types one arrives at the following
conditions:
\begin{itemize}
\item $\nabla F = 0$;
\item the Riemann curvature tensor is given by
  \begin{equation}
    \label{eq:riemd=11}
    \Riem(g) = \tfrac1{12} T^{[4]} + \tfrac1{36} (g \odot T^{[2]}) -
    \tfrac1{72} |F|^2 (g \odot g)~,
  \end{equation}
  where $\odot$ is the Kulkarni--Nomizu product (see, e.g.,
  \cite[§1.G, 1.110]{Besse}), and the tensors $T^{[2k]} \in C^\infty(M, S^2
  \Lambda^k T^*M)$, for $k=1,2$ are defined by
  \begin{equation*}
    \begin{aligned}[t]
      T^{[2]}(X,Y) &:= \left< \iota_X F, \iota_Y F\right>\\
      T^{[4]}(X,Y,W,Z) &:= \left<\iota_X \iota_Y F, \iota_W \iota_Z F
      \right>
    \end{aligned}
  \end{equation*}
  for all vector fields $X,Y,W,Z$; and
\item $F$ obeys the Plücker identity:
  \begin{equation*}
    \imath_X \imath_Y \imath_Z F \wedge F = 0~,
  \end{equation*}
  for all vector fields $X,Y,Z$.
\end{itemize}
The first two conditions imply that the Riemann tensor is parallel,
whence $(M,g)$ is locally symmetric, whence locally isometric to one
of the spaces in Theorem \ref{thm:CahenWallach}.  Every such space is
acted on transitively by a Lie group $G$ (the group of
\textbf{transvections}), whence if we fix a point in $M$ (the
\textbf{origin}) with isotropy $H$, $M$ is isomorphic to the space of
cosets $G/H$.  Let $\fg$ denote the Lie algebra of $G$ and $\fh$ the
Lie subalgebra corresponding to $H$.  Then $\fg$ admits a vector space
decomposition $\fg=\fh\oplus \fm$, where $\fm$ is isomorphic to the
tangent space of $M$ at the origin.  The Lie brackets are such that
\begin{equation*}
  [\fh,\fh] \subset \fh \qquad   [\fh,\fm] \subset \fm \qquad
  [\fm,\fm] \subset \fh~.
\end{equation*}
The metric $g$ on $M$ is determined by an $\fh$-invariant inner
product $B$ on $\fm$.

Since $F$ is parallel, it is $G$-invariant.  This means that it is
uniquely specified by its value at the origin, which defines an
$\fh$-invariant four-form on $\fm$.  For $F=0$, the right-hand side of
equation~\eqref{eq:riemd=11} vanishes, and hence $g$ is flat.  We will
therefore assume that $F\neq 0$.  The Plücker identity says that it is
then decomposable, whence it determines a four-dimensional vector
subspace $\fn \subset \fm$ as follows: if at the origin $F = \theta_1
\wedge \theta_2 \wedge \theta_3 \wedge \theta_4$, then $\fm$ is the
span of (the dual vectors to) the $\theta_i$.  Furthermore, because
$F$ is invariant, we have that $H$ leaves the space $\fn$ invariant,
whence $[\fh, \fn] \subset \fn$, which means that the holonomy group
of $M$ (which is isomorphic to $H$) acts reducibly.  In lorentzian
signature this does not imply that the space is locally isometric to a
product, since the metric may be degenerate when restricted to $\fn$.
Therefore we must distinguish between two cases, depending on whether
or not the restriction $B|_{\fn}$ of $B$ to $\fn$ is or is not
degenerate.

If $B|_{\fn}$ is non-degenerate, then it follows from the de~Rham--Wu
decomposition theorem \cite{Wu} that the space is locally isometric to
a product $N\times P$, with $N$ and $P$ locally symmetric spaces of
dimensions four and seven, respectively.  Explicitly, we can see this
as follows: there exists a $B$-orthogonal decomposition
$\fm=\fn\oplus\fp$, with $\fp := \fm^\perp$, where $[\fh,\fp] \subset
\fp$ because of the invariance of the inner product.  Let $\fg_N = \fh
\oplus \fn$ and $\fg_P = \fh \oplus \fp$.  They are clearly both Lie
subalgebras of $\fg$.  Let $G_N$ and $G_P$ denote the respective
(connected, simply-connected) Lie groups.  Then $N$ will be locally
isometric to $G_N/H$ and $P$ will be locally isometric to $G_P/H$, and
$M$ will be locally isometric to the product.  The metrics on $N$ and
$P$ are induced by the restrictions of $\fn$ and $\fp$ respectively of
the inner product $B$ on $\fn \oplus \fp$, denoted
\begin{equation}
  \begin{aligned}
    B_\fn&= B|_\fn
    \\
    B_\fp&=B|_\fp~.
  \end{aligned}
\end{equation}
We shall denote the metrics on $N$ and $P$ induced from the above
inner products by $h$ and $m$, respectively.

On the other hand if the restriction $B|_\fn$ is degenerate, so that
$\fn$ is a null four-dimensional subspace of $\fm$, the four-form $F$
is also null.  From Theorem~\ref{thm:CahenWallach} one sees (see,
e.g., \cite{FOPflux}) that the only lorentzian symmetric spaces
admitting parallel null forms are those which are locally isometric to
a product $M=\CW_d(A)\times Q_{11-d}$, where $\CW_d(A)$ is a
$d$-dimensional Cahen-Wallach space and $Q_{11-d}$ is an
$(11{-}d)$-dimensional riemannian symmetric space.

In summary, there are two separate cases to consider:
\begin{enumerate}
\item $(M,g)= (N_4 \times P_7, h \oplus m)$ (locally), where
  $(N,h)$ and $(P,m)$ are symmetric spaces and where $F$ is
  proportional to (the pull-back of) the volume form on $(N,h)$; or
\item $M = \CW_d(A) \times Q_{11-d}$ (locally) and $d\geq 3$, where
  $Q_{11-d}$ is a riemannian symmetric space.
\end{enumerate}

In \cite{FOPMax} these cases are analysed further, resulting in the
following theorem.

\begin{theorem}[FO-Papadopoulos \cite{FOPMax}]
  \label{th:d=11}
  Let $(M,g,F)$ be a maximally supersymmetric solution of
  eleven-dimensional supergravity.  Then it is locally isometric to
  one of the following:
  \begin{itemize}
  \item $\AdS_7(-7R) \times S^4(8R)$ and $F = \sqrt{6 R} \dvol(S^4)$,
    where $R>0$ is the constant scalar curvature of $M$;
  \item $\AdS_4(8R) \times S^7(-7R)$ and $F = \sqrt{-6 R}
    \dvol(\AdS_4)$, where $R<0$ is again the constant scalar curvature
    of $M$; or
  \item $\CW_{11}(A)$ with $A=-\frac{\mu^2}{36}\,
    \text{diag}(4,4,4,1,1,1,1,1,1)$ and\\
    $F = \mu\, dx^- \wedge dx^1\wedge dx^2 \wedge dx^3$.  One must
    distinguish between two cases:
    \begin{itemize}
    \item $\mu=0$: which recovers the flat space solution $\EE^{1,10}$
          with $F=0$; and
    \item $\mu\neq 0$: all these are isometric and describe a
      symmetric plane wave.
    \end{itemize}
  \end{itemize}
\end{theorem}

The first two solutions are the well-known Freund--Rubin backgrounds
\cite{FreundRubin} and \cite{AdS7S4}, whereas the plane wave was
originally discovered by Kowalski-Glikman \cite{KG} and rediscovered
subsequently in \cite{FOPflux}.  All of these solutions are locally
isometric to the intersection of two quadrics in $\EE^{11,2}$.
Moreover, as shown in \cite{ShortLimits,Limits} they are related by
``plane-wave limits'' \cite{PenrosePlaneWave,GuevenPlaneWave}.

\subsection{Ten-dimensional IIA supergravity}
\label{sec:maxIIA}

Type IIA supergravity \cite{GianiPerniciIIA, CampbellWestIIA,
  HuqNamazieIIA} is obtained by dimensional reduction from
eleven-dimensional supergravity.  From the discussion in
Section~\ref{sec:kk} and the fact that a $d=11$ background is
characterised by a metric $g$ and a $4$-form $F$, it follows that a
IIA supergravity background is characterised by a quintuplet
$(h,\phi,\Omega,G,H)$ where $h$ is a lorentzian metric on a
ten-dimensional spacetime $N$, $\phi$ is a real function on $N$,
$\Omega$ a closed $2$-form which is the curvature of a principal
$\Gamma$-bundle over $N$, $G$ a $4$-form on $N$ and $H$ a $3$-form on
$N$.  The PDEs satisfied by these fields are obtained by reducing
those in $d=11$ supergravity by the action of the group $\Gamma$.

It is a fundamental property of the Ka{\l}u\.za--Klein reduction, that
any IIA supergravity background can be lifted (or ``oxidised'') to a
background of eleven-dimensional supergravity possessing a
one-parameter group symmetries.  If the IIA supergravity solution
preserves some supersymmetry, its lift to eleven dimensions will
preserve at least the same amount of supersymmetry. This means that a
maximally supersymmetric solution of IIA supergravity will uplift to
one of the maximally supersymmetric solutions of eleven-dimensional
supergravity determined in the previous section.  Therefore the
determination of the maximally supersymmetric IIA backgrounds reduces
to classifying those dimensional reductions of the maximally
supersymmetric eleven-dimensional backgrounds which preserve all
supersymmetry.

As explained already in \cite{FOPflux}, the only such reductions are
the reductions of the flat eleven-dimensional background by a
translation subgroup of the Poincaré group.  In summary, one has

\begin{corollary}[FO-Papadopoulos \cite{FOPMax}]
  \label{th:IIA}
  Any maximally supersymmetric solution of type IIA supergravity is
  locally isometric to $\EE^{1,9}$ with zero fluxes and constant
  dilaton.
\end{corollary}

\subsection{Ten-dimensional IIB supergravity}
\label{sec:maxIIB}

A maximally supersymmetric background of IIB supergravity admits a
(real) $32$-dimensional space of Killing spinors.  Since this is the
(real) rank of the spinor bundle $\eS$ defined in
Section~\ref{sec:IIB}, it means that at any given point, there is a
basis for the spinor bundle consisting of Killing spinors.  These
spinors satisfy equation \eqref{eq:dilatinoIIB}, whence $c(B)$ and
$c(G)$ must vanish separately, which in turn imply the vanishing of
$G$ and $B$.  In particular, this has a consequence that $z$ and hence
$\tau$ are constant, whence the connection $A$ on $L_\tau$ also
vanishes.  Maximally supersymmetric backgrounds have the form
$(M,g,F)$ and are parametrised by the upper half-plane via the
constant parameter $\tau$.  Maximally supersymmetry now implies the
flatness of the connection $D$ defined by equation
\eqref{eq:gravitinoIIB} and which takes the simplified form
\begin{equation*}
  D_X \varepsilon = \nabla_X \varepsilon + \tfrac{i}4 c(F) c(X^\flat)
  \varepsilon~,
\end{equation*}
where $\nabla$ is the spin connection.  Notice that the equations of
motion now say that $F$ is closed.

Computing the curvature of this connection, and separating into types,
we arrive at the following conditions:
\begin{itemize}
\item $\nabla F = 0$;
\item the Riemann curvature tensor is given by
  \begin{equation}
    \label{eq:RiemIIB}
    R(X,Y,Z,W) = \left<\imath_X\imath_ZF,\imath_Y\imath_WF\right> - 
    \left<\imath_X\imath_WF,\imath_Y\imath_ZF\right>~.
  \end{equation}
  Since $F$ is parallel, this means that so is the Riemann tensor,
  whence $(M,g)$ is locally symmetric; and
\item $F$ obeys an identity reminiscent of both the Plücker and Jacobi
  identities:
  \begin{equation}
    \label{eq:plujac}
    \lambda(\imath_X\imath_Y\imath_Z F) F = 0 \qquad\text{for all
      vector fields $X,Y,Z$,}
  \end{equation}
  where $\lambda : \Omega^2(M) \to \End(\Lambda^5 T^*M)$ is the
  composition of the (metric-induced) isomorphism $\Omega^2(TM) \cong
  \fso(TM)$ between $2$-forms and skew-symmetric endomorphisms of the
  tangent bundle and the action of such endomorphisms on the
  $5$-forms.
\end{itemize}

It was proved in \cite{FOPPluecker} that equation~\eqref{eq:plujac}
implies that $F = G + \star G$, where $G = \theta_1 \wedge \theta_2
\wedge \theta_3 \wedge \theta_4 \wedge \theta_5$ is a parallel
decomposable form.

The ensuing analysis follows closely the case of eleven-dimensional
supergravity and will not be repeated here.  We must distinguish
between two cases, depending on whether or not the five-form $G$ is
null.  First suppose that $G$ (and hence $F$) is not null.  Then the
five-form $G$ induces a local decomposition of $(M,g)$ into a product
$N_5\times P_5$ of two five-dimensional symmetric spaces $(N,h)$ and
$(P,m)$, where $G \propto \dvol(N)$ and hence $\star G \propto
\dvol(P)$.  Since $(M,g)$ is lorentzian, one of the spaces $(N,h)$ and
$(P,m)$ is lorentzian and the other riemannian.  By interchanging $G$
with $\star G$ if necessary, we can assume that $G$ has positive norm
and hence that $N$ is riemannian.

In summary, there are two separate cases to consider:
\begin{enumerate}
\item $(M,g)= (N_5 \times P_5, h \oplus m)$ (locally), where
  $(N,h)$ and $(P,m)$ are symmetric spaces and where $F=G + \star G$
  and $G$ is proportional to (the pull-back of) the volume form on
  $(N,h)$; or
\item $M = \CW_d(A) \times Q_{10-d}$ (locally) and $d\geq 3$, where
  $Q_{10-d}$ is a riemannian symmetric space.
\end{enumerate}

In \cite{FOPMax} these cases are analysed further, resulting in the
following theorem.

\begin{theorem}[FO-Papadopoulos \cite{FOPMax}]
  \label{th:IIB}
  Let $(M,g,F_5^+,...)$ be a maximally supersymmetric solution of
  ten-dimensional type IIB supergravity.  Then it has constant
  axi-dilaton (normalised so that $z=0$ in the formulas below), all fluxes 
  vanish except for the one corresponding to the self-dual five-form,
  and is locally isometric to one of the following:
  \begin{itemize}
  \item $\AdS_5(-R) \times S^5(R)$ and $F = 2 \sqrt{\frac{R}5}
    \left(\dvol(\AdS_5) + \dvol(S^5)\right)$, where $\pm R$ are the
    scalar curvatures of $\AdS_5$ and $S^5$, respectively; or
  \item $\CW_{10}(A)$ with $A=-\mu^2 \1$ and $F = \half \mu\, dx^-
    \wedge (dx^1\wedge dx^2 \wedge dx^3 \wedge dx^4 + dx^5\wedge dx^6
    \wedge dx^7 \wedge dx^8)$.  One must distinguish between two
    cases:
    \begin{itemize}
    \item $\mu=0$: which yields the flat space solution $\EE^{1,9}$
      with zero fluxes; and
    \item $\mu\neq 0$: all these are isometric and describe a
      symmetric plane wave.
    \end{itemize}
  \end{itemize}
\end{theorem}

The first solution is the well-known Freund--Rubin background
mentioned originally in \cite{SchwarzIIB}.  The plane wave solution
was discovered in \cite{NewIIB}.  As in eleven-dimensional
supergravity, the solutions above are locally isometric to the
intersection of two quadrics in $\EE^{10,2}$ and as shown in
\cite{ShortLimits,Limits} they are related by plane-wave limits.

\subsection{Six-dimensional $(2,0)$ and $(1,0)$ supergravities}
\label{sec:maxsixd}

In this case, maximal supersymmetry implies the flatness of the
supersymmetric connection $D$ in \eqref{eq:conn10} which, as explained
in Section~\ref{sec:sixdim}, is induced from a metric connection with
closedd torsion $3$-form $H$.  In Section~\ref{sec:flat} we showed
that $(M,g)$ is locally isometric to a six-dimensional Lie group with
a bi-invariant lorentzian metric.  The only extra condition is that
$H$, the canonical bi-invariant $3$-form associated to such a Lie
group, should be self-dual.

We therefore look for Lie algebras with invariant lorentzian scalar
products relative to which the canonical invariant $3$-form is
anti-self dual.  As explained in Section~\ref{sec:lie}, such Lie
algebra is a direct sum of indecomposables.  Furthermore, if the Lie
algebra is indecomposable then it must be the double extension of an
abelian Lie algebra by a one-dimensional Lie algebra and hence
solvable (see, e.g., \cite{MedinaRevoy}).

These considerations make possible the following enumeration of
six-dimensional lorentzian Lie algebras:
\begin{enumerate}
\item $\EE^{1,5}$
\item $\EE^{1,2} \oplus \fso(3)$
\item $\EE^3 \oplus \fso(1,2)$
\item $\fso(1,2) \oplus \fso(3)$
\item $\fd(\EE^4,\RR)$
\end{enumerate}
where the last case actually corresponds to a family of Lie algebras,
depending on the action of $\RR$ on $\EE^4$, which is given by a
homomorphism $\RR \to \fso(4)$.

Imposing the condition of anti-selfduality trivially discards cases
(2) and (3) above.  Case (1) is the abelian Lie algebra with Minkowski
metric.  The remaining two cases were investigated in
\cite{CFOSchiral} (see also \cite{GMR}) in detail and we review this
below.

\subsubsection{A six-dimensional Cahen--Wallach space}

Let $e_i$, $i=1,2,3,4$, be an orthonormal basis for $\EE^4$, and let
$e_- \in \RR$ and $e_+\in\RR^*$, so that together they span
$\fd(\EE^4,\RR)$.  The action of $\RR$ on $\EE^4$ defines a map $\rho:
\RR \to \Lambda^2 \EE^4$, which can be brought to the form $\rho(e_-) =
\alpha e_1 \wedge e_2 + \beta e_3 \wedge e_4$ via an orthogonal change
of basis in $\EE^4$ which moreover preserves the orientation.  The Lie
brackets of $\fd(\EE^4,\RR)$ are given by
\begin{equation*}
  \begin{aligned}[m]
    [e_-,e_1] &= \alpha e_2\\
    [e_-,e_2] &= -\alpha e_1\\
    [e_1,e_2] &= \alpha e_+
  \end{aligned}\qquad\qquad
  \begin{aligned}[m]
    [e_-,e_3] &= \beta e_4\\
    [e_-,e_4] &= -\beta e_3\\
    [e_3,e_4] &= \beta e_+
  \end{aligned}~,
\end{equation*}
and the scalar product is given (up to scale) by
\begin{equation*}
  \left<e_-,e_-\right> = b \qquad \left<e_+,e_-\right> = 1 \qquad
  \left<e_i,e_j\right> = \delta_{ij}~.
\end{equation*}
The first thing we notice is that we can set $b=0$ without loss of
generality by the automorphism fixing all $e_i,e_+$ and mapping $e_-
\mapsto e_- - \half b e_+$.  We will assume that this has been done
and that $\left<e_-,e_-\right>=0$.  A straightforward calculation
shows that the three-form $H$ is anti-selfdual if and only if $\beta =
\alpha$.  Let us put $\beta=\alpha$ from now on.  We must distinguish
between two cases: if $\alpha = 0$, then the resulting algebra is
abelian and is precisely $\EE^{1,5}$.  On the other hand if
$\alpha\neq 0$, then rescaling $e_\pm \mapsto \alpha^{\pm 1} e_\pm$ we
can effectively set $\alpha = 1$ without changing the scalar product.
Finally we notice that a constant rescaling of the scalar product can
be undone by an automorphism of the algebra.  As a result we have two
cases: $\EE^{1,5}$ (obtained from $\alpha =0$) and the algebra
\begin{equation}
  \label{eq:nw6}
  \begin{aligned}[m]
    [e_-,e_1] &= e_2\\
    [e_-,e_2] &= - e_1\\
    [e_1,e_2] &= e_+
  \end{aligned}\qquad\qquad
  \begin{aligned}[m]
    [e_-,e_3] &= e_4\\
    [e_-,e_4] &= - e_3\\
    [e_3,e_4] &= e_+
  \end{aligned}~,
\end{equation}
with scalar product given by
\begin{equation}
  \label{eq:nw6m}
  \left<e_+,e_-\right> = 1 \qquad\text{and}\qquad
  \left<e_i,e_j\right> = \delta_{ij}~.
\end{equation}
There is a unique simply-connected Lie group with the above Lie
algebra which inherits a bi-invariant lorentzian metric.  This Lie
group is a six-dimensional analogue of the Nappi--Witten group
\cite{NW}, which is based on the double extension $\fd(\EE^2,\RR)$
\cite{FSsug}.  This was denoted $\NW_6$ in \cite{FSPL}, where one
can find a derivation of the metric on this six-dimensional group.
The supergravity solution was discovered by Meessen \cite{Meessen} who
called it KG6 by analogy with the maximally supersymmetric plane wave
of eleven-dimensional supergravity discovered by Kowalski-Glikman
\cite{KG} and rediscovered in \cite{FOPflux}.

The metric is easy to write down once we choose a parametrisation for
the group.  The calculation is routine (see, for example, \cite{FSPL})
and the result is
\begin{equation}
  \label{eq:CWmetric2}
  g = 2 dx^+ dx^- - \tfrac14 \sum_i (x^i)^2 (dx^-)^2 + \sum_i
  (dx^i)^2~.
\end{equation}
In these coordinates the three-form $H$ is given by
\begin{equation*}
  H = \tfrac23 dx^- \wedge (dx^1 \wedge dx^2 + dx^3 + dx^4)~.
\end{equation*}

\subsubsection{The Freund--Rubin backgrounds}

Finally we discuss case (4), with Lie algebra $\fso(1,2) \oplus
\fso(3)$.  Let $e_0,e_1,e_2$ be a pseudo-orthonormal basis for
$\fso(1,2)$.  The Lie brackets are given by
\begin{equation*}
  [e_0,e_1] = -e_2\qquad   [e_0,e_2] = e_1\qquad   [e_1,e_2] =
  e_0~.
\end{equation*}
Similarly let $e_3,e_4,e_5$ denote an orthonormal basis for
$\fso(3)$, with Lie brackets
\begin{equation*}
  [e_5,e_3] = -e_4\qquad   [e_5,e_4] = e_3\qquad   [e_3,e_4] =
  -e_5~.
\end{equation*}
The most general invariant lorentzian scalar product on
$\fso(1,2)\oplus\fso(3)$ is labelled by two positive numbers $\alpha$
and $\beta$ and is given by
\begin{equation*}
  \bordermatrix{& e_0 & e_1 & e_2 & e_3 & e_4 & e_5 \cr
  e_0 &  -\alpha & 0  &  0  & 0 & 0 & 0 \cr
  e_1 & 0 & \alpha & 0  & 0 & 0 & 0 \cr
  e_2 & 0 & 0 & \alpha & 0 & 0 & 0 \cr
  e_3 & 0 & 0 & 0 & \beta & 0 & 0 \cr
  e_4 & 0 & 0 & 0 & 0 & \beta &  0 \cr
  e_5 & 0 & 0 & 0 & 0 & 0 & \beta \cr}~.  
\end{equation*}
Anti-selfduality of the canonical three-form implies that $\beta =
\alpha$.  There is a unique simply-connected Lie group with Lie
algebra $\fso(1,2) \oplus \fso(3)$, namely $\widetilde{\SL(2,\RR)}
\times \SU(2)$, where $\widetilde{\SL(2,\RR)}$ denotes the universal
covering group of $\SL(2,\RR)$.  This group inherits a one-parameter
family of bi-invariant metrics.  This solution is none other than the
standard Freund--Rubin solution $\AdS_3 \times S^3$, with equal radii
of curvature, where strictly speaking we should take the universal
covering space of $\AdS_3$.

In summary, the following are the possible maximally supersymmetric
backgrounds of $(1,0)$ supergravity, and of $(2,0)$ supergravity up to
the action of the R-symmetry group.  First of all we have a
one-parameter family of Freund-Rubin backgrounds locally isometric to
$\AdS_3 \times S^3$, with equal radii of curvature.  The anti-selfdual
three-form $H$ is then proportional to the difference of the volume
forms of the two spaces.  Then we have a six-dimensional analogue
$\NW_6$ of the Nappi--Witten group, locally isometric to a
Cahen--Wallach symmetric space.  Finally there is flat Minkowski
spaceteime $\EE^{1,5}$.  These backgrounds are related by Penrose
limits which can be interpreted in this case as group contractions.
The details appear in \cite{FSPL}.

\subsection{Five-dimensional $N=2$ supergravity}
\label{sec:maxN=2}

In this section we will review the dimensional reductions of the
six-dimensional backgrounds just found.  Dimensional reduction usually
breaks some supersymmetry: in the ten- and eleven-dimensional
supergravity theories, only the flat background remains maximally
supersymmetric after dimensional reduction and then only by a
translation.  However for the six-dimensional backgrounds the
situation is different.  Indeed, in \cite{LMO8} it was shown that the
thereto known maximally supersymmetric backgrounds with eight
supercharges in six, five and four dimensions are related by
dimensional reduction and oxidation.  As we will see presently, this
perhaps surprising phenomenon stems from the fact that the
six-dimensional backgrounds are parallelised Lie groups.  Our results
will also give an \emph{a priori} explanation to the empirical fact
that these backgrounds are homogeneous \cite{ALO}.

We now explain the technical result which underlies this result.  Let
$D$ be a metric connection with torsion $T$.  We observe that if a
vector field $\xi$ is $D$-parallel then it is Killing.  Now let $\psi$
be a Killing spinor; that is, $D\psi = 0$. Then the Lie derivative of
$\psi$ along $\xi$ is well-defined (see, for example,
\cite{JMFKilling}) and, furthermore, it vanishes identically.
Moreover, if $\eL_\xi \psi = 0$ for \emph{all} Killing spinors then
$D\xi = 0$.

For a parallelised Lie group $G$, the $D$-parallel vectors are either
the left- or right-invariant vector fields, depending on the choice of
parallelising connection.  For definiteness, we will choose the
connection whose parallel sections are the left-invariant vector
fields.  Left-invariant vector fields generate right translations and
are in one-to-one correspondence with elements of the Lie algebra
$\fg$.  Therefore every left-invariant vector field $\xi$ determines a
one-parameter subgroup $K$, say, of $G$ and the orbits of such a
vector field in $G$ are the right $K$-cosets.  The dimensional
reduction along this vector field is smooth and diffeomorphic to the
space of cosets $G/K$.  We will be interested in subgroups $K$ such
that $G/K$ is a five-dimensional lorentzian spacetime, which requires
that the right $K$-cosets are spacelike.  In other words, we require
that the Killing vector $\xi$ be spacelike.  Bi-invariance of the
metric guarantees that this is the case provided that the Lie algebra
element $\xi(e)\in \fg$ is spacelike relative to the invariant scalar
product.  Further notice that a constant rescaling of $\xi$ does not
change its causal property nor the subgroup $K$ it generates: it is
simply reparameterised.  Therefore, in order to classify all possible
reductions we need to classify all spacelike elements of $\fg$ up to
scale.  Moreover elements of $\fg$ which are related by isometric
automorphisms (e.g., which are in the same adjoint orbit of $G$) give
rise to isometric quotients.  Thus, to summarise, we want to classify
spacelike elements of $\fg$ up to scale and up to automorphisms.

As discussed in Section~\ref{sec:kk}, the reduction of the
six-dimensional metric to five dimensions gives rise to a metric $h$,
a dilaton $\phi$ and a curvature 2-form $F$.  The dilaton $\phi$ is a
logarithmic measure of the fibre metric $\|\xi_X\|$ which in our case
is constant, and $F = d\alpha$ (omitting pullbacks).  We can give an
explicit formula for $\Omega$ using the Maurer--Cartan structure
equations.  Indeed,
\begin{equation}
  \label{eq:F}
  F = d\alpha = \left<X,d\theta\right> = -\half
  \left<X,[\theta,\theta]\right>~.
\end{equation}
In terms of this data, the metric on the $G$ is given by the usual
Ka{\l}u\.za--Klein ansatz
\begin{equation*}
  ds^2 = h + \alpha^2~,
\end{equation*}
where we have set the dilaton to zero in agreement with the choice of
normalisation for $\xi_X$.  More explicitly the metric on the
five-dimensional quotient is given by
\begin{equation*}
  h = \left<\theta,\theta\right> - \left<X,\theta\right>^2~.
\end{equation*}

To reduce the anti-selfdual three-form $H$ we first decompose it as
\begin{equation*}
  H = G_3 + \alpha \wedge G_2~,
\end{equation*}
where $G_2 = \iota_{\xi_X} H$ and $G_3$ are basic.  Because $dH=0$ it
follows that $dG_2=0$ and that $dG_3 + F \wedge G_2=0$ where $F =
d\alpha$ was defined above.  Finally because $H$ is anti-selfdual, it
follows that $G_3$ and $G_2$ are related by Hodge duality in five
dimensions: $G_3 = \star_h G_2$.  In other words, we have that
\begin{equation*}
  H = \star_h G_2 + \alpha \wedge G_2~,
\end{equation*}
where $dG_2 = 0$ and $d\star_h G_2 = -F \wedge G_2$.

In fact, in this case we have $F=G_2$.  Indeed, using that $H = -\tfrac16
\left<\theta,[\theta,\theta]\right>$, we compute
\begin{equation*}
   G_2 = \iota_{\xi_X} H = -\half \left<X,[\theta,\theta]\right>~,
\end{equation*}
which agrees with the expression for $F$ derived in \eqref{eq:F}.

In summary, for the reductions under consideration, we obtain a
maximally supersymmetric background of the minimal $N{=}2$
supergravity with bosonic fields $(h,F)$ given by the reduction of
$(g,H)$ where $F=d\alpha$, $h = g - \alpha^2$ and $H = \star_h F +
\alpha \wedge F$.

The different reductions were classified in \cite{CFOSchiral}, to
where we send the reader for details, hence obtaining all the
maximally supersymmetric backgrounds of the minimal $N{=}2$
supergravity and thus completing the classification of supersymmetric
backgrounds in \cite{GGHPR}.  Among the maximally supersymmetric
backgrounds one finds the near-horizon geometries \cite{GMT1} of the
rotating black holes of \cite{BMPV,KRW}, the symmetric plane-wave of
\cite{Meessen} and the Gödel-like background discovered in
\cite{GGHPR}.

\section{Parallelisable type II backgrounds}
\label{sec:parallel}

In this section we will present a classification of parallelisable
type II backgrounds,  by which we mean backgrounds of both type IIA
and type IIB supergravity.  Since these theories contain different
dynamical degrees of freedom, common backgrounds are necessarily very
special.

\begin{definition}
  A type II supergravity
  \idx{background}{supergravity!$d{=}10$ type~II!background} consists of a
  ten-dimensional lorentzian spin manifold $(M,g)$ together with a
  closed $3$-form $H$ and a smooth function $\phi: M \to \RR$ subject
  to the equations of motion obtained by varying the (formal) action
  functional
  \begin{equation}
    \label{eq:typeIIaction}
    \int_M e^{-2\phi} \left( R + 4 |d\phi|^2 - \half |H|^2 \right)
    \dvol_g ~,
  \end{equation}
  where $R$ and $\dvol_g$ are the scalar curvature and the volume form
  associated to $g$.  
\end{definition}

We are interested in
\idx{parallelisable}{supergravity!$d{=}10$ type~II!parallelisable background}
backgrounds, for which the metric connection $D$ with torsion $3$-form
$H$ is flat.  In that case, the equations of motion simplify to the
following three conditions:
\begin{equation}
  \label{eq:eomspara}
  \begin{aligned}[c]
    \nabla d \phi &= 0\\
    d\phi \wedge \star H &= 0\\
    |d\phi|^2 - \tfrac14 |H|^2 &=0~.
  \end{aligned}
\end{equation}

To discuss supersymmetry, we need to distinguish whether we are in
type IIA or type IIB supergravity, since the spinor bundles are
different.  Let $S_\pm$ denote the real $16$-dimensional half-spin
representations of $\Spin(1,9)$ and let $S_A = S_+ \oplus S_-$ and
$S_B = S_+ \oplus S_+$.  Let $\eS_A$ and $\eS_B$ denote the spinor
bundles on $M$ associated to $S_A$ and $S_B$, respectively.  We will
let $\eS$ denote either $\eS_A$ or $\eS_B$, depending on which type II
theory we are considering.

\begin{definition}
  A type II background is
  \idx{supersymmetric}{supergravity!$d{=}10$ type~II!supersymmetric background}
  if there are nonzero sections $\varepsilon$ of $\eS$ satisfying the
  two conditions:
  \begin{equation*}
    D \varepsilon = 0 \qquad\text{and}\qquad c\left(d\phi + \half
      H\right) \varepsilon = 0~,
  \end{equation*}
  where $c: \Omega(M) \to \Cl(TM) \to \End(\eS)$ is the Clifford
  action of forms on spinors.
\end{definition}

The supersymmetric parallelisable type II backgrounds were classified
in \cite{JMFPara,KYPara} and revisited in the context of heterotic
supergravity in \cite{FKYHeterotic}, whose treatment we follow.

\subsection{Ten-dimensional parallelisable geometries}
\label{sec:parallelisms}

As explained in Section \ref{sec:parallelisable}, it is possible to
list all the simply-connected parallelisable lorentzian manifolds in
any dimension.  The ingredients out of which we can make them are
given in Table~\ref{tab:ingredients}, whose last column follows from
equation \eqref{eq:eomspara}.

Indeed, in the case of a Lie group, that is, when $dH=0$, equation
\eqref{eq:eomspara} says that $d\phi$ must be central, when thought of as
an element in the Lie algebra. Since $\AdS_3$, $S^3$ and $\SU(3)$ are
simple, their Lie algebras have no centre, whence $d\phi=0$.  In the
case of an abelian group there are no conditions, and in the case of
$\CW_{2n}(A)$, the Lie algebra has a one-dimensional centre
corresponding to $\d_+$, whose dual one-form is $dx^-$.  This means
that $d\phi$ must be proportional to $dx^-$, whence $\phi$ can only
depend on $x^-$.  Finally for $S^7$, the equation of motion $\star H
\wedge d\phi = 0$ implies that $d\phi=0$.  To see this, notice that
the parallelised $S^7$ possesses a nearly parallel $G_2$ structure and
the differential forms decompose into irreducible types under $G_2$.
For example, the one-forms corresponding to the irreducible
seven-dimensional irreducible representation $\fm$ of $G_2$ coming
from the embedding $G_2 \subset \SO(7)$, whereas the two-forms
decompose into $\fg_2 \oplus \fm$, where $\fg_2$ is the adjoint
representation which is irreducible since $G_2$ is simple.  Now, $H$
and $\star H$ both are $G_2$-invariant and hence the map
$\Omega^1(S^7) \to \Omega^2(S^7)$ defined by $\theta \mapsto \star
(\star H \wedge\theta)$ is $G_2$-equivariant.  Since it is not
identically zero, it must be an isomorphism onto its image.  Hence if
$\star H \wedge d\phi = 0$, then also in this case $d\phi = 0$.

\begin{table}[h!]
  \centering
  \setlength{\extrarowheight}{3pt}
  \renewcommand{\arraystretch}{1.3}
  \begin{small}
    \begin{tabular}{|>{$}l<{$}|>{$}l<{$}|>{$}l<{$}|}\hline
      \multicolumn{1}{|c|}{Space} & \multicolumn{1}{c|}{Torsion} &
      \multicolumn{1}{c|}{Dilaton} \\
      \hline\hline
      \AdS_3 & dH=0 \quad |H|^2 < 0 & \text{constant}\\
      \EE^{1,0} & H=0 & \text{unconstrained}\\
      \EE^{0,1} & H=0 & \text{unconstrained}\\
      S^3 & dH=0 \quad |H|^2 > 0 &  \text{constant}\\
      S^7 & dH\neq 0 \quad |H|^2 > 0 & \text{constant}\\
      \SU(3) & dH= 0 \quad |H|^2 > 0 & \text{constant}\\
      \CW_{2n}(A) & dH=0 \quad |H|^2 = 0 & \phi(x^-)\\ \hline
    \end{tabular}
  \end{small}
  \vspace{8pt}
  \caption{Elementary parallelisable (lorentzian or riemannian) geometries}
  \label{tab:ingredients}
\end{table}

It is now a simple matter to put these ingredients together to make up
all possible ten-dimensional combinations with lorentzian signature.
Doing so, we arrive at Table~\ref{tab:geometries} (see also
\cite{JMFPara}, where the entry corresponding to $\EE^{1,0} \times S^3
\times S^3 \times S^3$ had been omitted inadvertently and where the
entries with $S^7$ had also been omitted due to the fact that in type
II string theory $dH=0$).

\begin{table}[h!]
  \centering
  \setlength{\extrarowheight}{3pt}
  \renewcommand{\arraystretch}{1.3}
  \begin{small}
    \begin{tabular}{|>{$}l<{$}|>{$}l<{$}|}\hline
      \multicolumn{1}{|c|}{Spacetime} & \multicolumn{1}{c|}{Spacetime}\\
      \hline\hline
      \AdS_3 \times S^7 & \AdS_3 \times S^3 \times S^3 \times \EE\\
      \AdS_3 \times S^3 \times \EE^4 & \AdS_3 \times \EE^7\\
      \EE^{1,0} \times S^3 \times S^3 \times S^3 & \EE^{1,1} \times \SU(3)\\
      \EE^{1,2} \times S^7 & \EE^{1,3} \times S^3 \times S^3\\
      \EE^{1,6} \times S^3 & \EE^{1,9}\\
      \CW_{10}(A) & \CW_8(A) \times \EE^2\\
      \CW_6(A) \times S^3 \times \EE & \CW_6(A) \times \EE^4\\
      \CW_4(A) \times S^3 \times S^3 & \CW_4(A) \times
      S^3 \times \EE^3\\
      \CW_4(A) \times \EE^6 & \\\hline
    \end{tabular}
  \end{small}
  \vspace{8pt}
  \caption{Ten-dimensional simply-connected parallelisable spacetimes}
  \label{tab:geometries}
\end{table}

\subsection{Type II backgrounds}
\label{sec:linear}

First of all we notice that $S^7$ cannot appear because $dH=0$.
Therefore the allowed backgrounds follow \emph{mutatis mutandis} from
the analysis of \cite{JMFPara,KYPara}.  We start by listing the
possible backgrounds and then counting the amount of supersymmetry
that each preserves.  The results are summarised in
Table~\ref{tab:linear} and Table~\ref{tab:summary}, which also
contains the analysis of the supersymmetry preserved by the background.

\subsubsection*{$\AdS_3 \times S^3 \times S^3 \times \EE$}

Here $d\phi$ can only have nonzero components along the flat
direction, which is spacelike, whence $|d\phi|^2 \geq 0$.  Equation
\eqref{eq:eomspara} says that $|H|^2 \geq 0$, so that if we call $R_0$,
$R_1$ and $R_2$ the radii of curvature of $\AdS_3$ and of the two
$3$-spheres, respectively, then
\begin{equation*}
  \frac{1}{R_1^2}  +  \frac{1}{R_2^2} \geq \frac{1}{R_0^2}~.
\end{equation*}
This bound is saturated if and only if the dilaton is constant.

\subsubsection*{$\AdS_3 \times S^3 \times \EE^4$}

This is the limit $R_2 \to \infty$ of the above case.

\subsubsection*{$\AdS_3 \times \EE^7$}

This would be the limit $R_1 \to \infty$ of the above case, but then
the inequality $R_0^{-2} \leq 0$ cannot be satisfied.  Hence this
geometry is not a background (with or without supersymmetry).

\subsubsection*{$\EE^{1,9}$}

In this case $H=0$, so $|d\phi|^2 =0$.  So we can take a linear
dilaton along a null direction: $\phi = a + b x^-$, for some constants
$a,b$ say.

\subsubsection*{$\EE^{1,0} \times S^3 \times S^3 \times S^3$}

The dilaton can only depend on the flat coordinate, which is timelike,
so $|d\phi|^2 \leq 0$.  However $|H|^2 > 0$, whence this geometry
is never a background (with or without supersymmetry).

\subsubsection*{$\EE^{1,1} \times \SU(3)$}

Here $|H|^2 > 0$, and $d\phi$ can have components along
$\EE^{1,1}$.  Letting $(x^0,x^1)$ be flat coordinates for $\EE^{1,1}$,
we can take $\phi = a + \half |H| x^1$, for some constant $a$,
without loss of generality.

\subsubsection*{$\EE^{1,3} \times S^3 \times S^3$}

Here $|H|^2 > 0$ and $d\phi$ can have components along
$\EE^{1,3}$ \cite{Khuri}.  With $(x^0,x^1,x^2,x^3)$ being flat
coordinates for $\EE^{1,3}$, we take $\phi = a + \half |H| x^1$, for
some constant $a$.

\subsubsection*{$\EE^{1,6} \times S^3$}

This is the limit $R_2 \to \infty$ of the above case, where $R_2$ is
the radius of curvature of one of the spheres
\cite{DuffLu5Brane,NS5Brane}.

\subsubsection*{$\CW_{2n}(A) \times \EE^{10-2n}$, $n=2,3,4,5$}

In these cases $|H|^2 = 0$ and hence $|d\phi|^2 = 0$, so that it
cannot have components along the flat directions (if any).  This means
$\phi = a + b x^-$, for constants $a,b$.

\subsubsection*{$\CW_{4}(A) \times S^3 \times S^3$}

Here $|d\phi|^2 = 0$, whereas $|H|^2 >0$, hence there are no
backgrounds with this geometry.

\subsubsection*{$\CW_{2n}(A) \times S^3 \times \EE^{7-2n}$,
  $n=2,3$}

Here $|H|^2 > 0$, whence $|d\phi|^2 > 0$.  This means that we can
take $\phi = a + b x^- + \half |H| y$, where $y$ is any flat
coordinate in $\EE^{7-2n}$ and $a,b$ are constants.

\begin{table}[h!]
  \centering
  \setlength{\extrarowheight}{3pt}
  \renewcommand{\arraystretch}{1.3}
  \begin{small}
    \begin{tabular}{|>{$}l<{$}|>{$}l<{$}|}\hline
      \multicolumn{1}{|c|}{Geometry} &
      \multicolumn{1}{c|}{Dilaton} \\
      \hline\hline
      \AdS_3 \times S^3 \times S^3 \times \EE & \phi = a + \half |H| y\\
      \AdS_3 \times S^3 \times \EE^4 & \phi = a + \half |H| y\\
      \EE^{1,1} \times \SU(3) &  \phi = a + \half |H| y\\
      \EE^{1,3} \times S^3 \times S^3 & \phi = a + \half |H| y\\
      \EE^{1,6} \times S^3 & \phi = a + \half |H| y\\
      \EE^{1,9} & \phi = a + b x^-\\
      \CW_{10}(A) & \phi = a + b x^-\\
      \CW_8(A) \times \EE^2 & \phi = a + b x^-\\
      \CW_6(A) \times S^3 \times \EE & \phi = a + b x^- + \half
      |H| y\\
      \CW_6(A) \times \EE^4 & \phi = a + b x^-\\
      \CW_4(A) \times S^3 \times \EE^3 & \phi = a + b x^- +
      \half |H| y\\
      \CW_4(A) \times \EE^6 & \phi = a + b x^-\\
      \hline
    \end{tabular}
  \end{small}
  \vspace{8pt}
  \caption{Parallelisable backgrounds with a linear dilaton.  The
    notation is such that $y$ is a spacelike flat coordinate.}
  \label{tab:linear}
\end{table}

\begin{table}[h!]
  \centering
  \setlength{\extrarowheight}{3pt}
  \renewcommand{\arraystretch}{1.3}
  \begin{small}
    \begin{tabular}{|>{$}l<{$}|>{$}c<{$}|>{$}c<{$}|}\hline
      \multicolumn{1}{|c|}{Parallelisable} &
      \multicolumn{2}{c|}{Supersymmetries with dilaton being}\\
      \multicolumn{1}{|c|}{geometry} &
      \multicolumn{1}{c|}{constant} &
      \multicolumn{1}{c|}{nonconstant} \\
      \hline\hline
      \AdS_3 \times S^3 \times S^3 \times \EE &  16 & 16 \\
      \AdS_3 \times S^3 \times \EE^4 & 16 & 16  \\
      \EE^{1,1} \times \SU(3) & \times & 16 \\
      \EE^{1,3} \times S^3 \times S^3 & \times & 16\\
      \EE^{1,6} \times S^3 & \times & 16 \\
      \EE^{1,9} & 32 & 16 \\
      \CW_{10}(A) & 16, 18(A), 20, 22(A), 24(B), 28(B) & 16 \\
      \CW_8(A) \times \EE^2& 16, 20 & 16\\
      \CW_6(A) \times S^3 \times \EE & \times & 16\\
      \CW_6(A) \times \EE^4 & 16, 24 & 16 \\
      \CW_4(A) \times S^3 \times \EE^3 & \times & 16 \\
      \CW_4(A) \times \EE^6 & 16 & 16\\
      \hline
    \end{tabular}
  \end{small}
  \vspace{8pt}
  \caption{Supersymmetric parallelisable backgrounds, with (A) or (B)
    indicating IIA or IIB.}
  \label{tab:summary}
\end{table}

\bibliographystyle{utphys}
\bibliography{AdS3,AdS,ESYM,Sugra,Geometry}


\begin{theindex}

  \item Lie algebra
    \subitem double extension, 10
    \subitem indecomposable, 10
    \subitem lorentzian, 9
  \item lorentzian symmetric space
    \subitem anti\nobreakspace  {}de\nobreakspace  {}Sitter, 3
    \subitem Cahen--Wallach, 3
    \subitem de\nobreakspace  {}Sitter, 3
    \subitem isometric embedding, 5
    \subitem parallelisable, 6

  \indexspace

  \item supergravity
    \subitem $d{=}10$ IIB
      \subsubitem background, 16
      \subsubitem duality group, 17
      \subsubitem Killing spinor, 17
    \subitem $d{=}10$ type\nobreakspace  {}II
      \subsubitem background, 29
      \subsubitem parallelisable background, 29
      \subsubitem supersymmetric background, 29
    \subitem $d{=}11$
      \subsubitem background, 14
      \subsubitem Killing spinor, 14
      \subsubitem supersymmetric background, 14
    \subitem $d{=}6$
      \subsubitem $(1,0)$ background, 18
      \subsubitem $(1,0)$ supersymmetric background,\nobreakspace 18
      \subsubitem $(2,0)$ background, 18
      \subsubitem $(2,0)$ supersymmetric background,\nobreakspace 18
      \subsubitem R-symmetry groups, 18

\end{theindex}

\end{document}